\documentclass[11pt]{article}
\usepackage{amscd}
\usepackage{}
\usepackage{mathrsfs}
\usepackage{amsfonts}
\usepackage{amsfonts,amssymb,mathrsfs}
\usepackage{amsmath,amscd}
\usepackage[dvips]{graphicx}
\usepackage[all]{xy}
\usepackage{graphicx}
\usepackage{fancyhdr}
\usepackage{mathptm,pslatex}
\usepackage{amsthm}

\begin{document}
\sloppy
\newcommand{\dickebox}{{\vrule height5pt width5pt depth0pt}}
\newtheorem{Def}{Definition}[section]
\newtheorem{Bsp}{Example}[section]
\newtheorem{Prop}[Def]{Proposition}
\newtheorem{Theo}[Def]{Theorem}
\newtheorem{Lem}[Def]{Lemma}
\newtheorem{Koro}[Def]{Corollary}
\newcommand{\overpr}{$\hfill\square$}
\newcommand{\pd}{{\rm proj.dim\, }}
\newcommand{\id}{{\rm inj.dim\, }}
\newcommand{\fd}{{\rm fin.dim}\,}
\newcommand{\add}{{\rm add \, }}
\newcommand{\Hom}{{\rm Hom \, }}
\newcommand{\gldim}{{\rm gl.dim}\,}
\newcommand{\End}{{\rm End \, }}
\newcommand{\Ext}{{\rm Ext}}
\newcommand{\D}{\rm D \,}
\newcommand{\Coker}{{\rm Coker}\,\,}
\newcommand{\cpx}[1]{#1^{\bullet}}
\newcommand{\Dz}[1]{{\rm D}^+(#1)}
\newcommand{\Df}[1]{{\rm D}^-(#1)}
\newcommand{\DFz}[1]{{\rm D_{F}}^+(#1)}
\newcommand{\DFf}[1]{{\rm D_{F}}^-(#1)}
\newcommand{\Db}[1]{{\rm D^b}(#1)}
\newcommand{\DFb}[1]{{\rm D_{F}^b}(#1)}
\newcommand{\DF}[1]{{\rm D_{F}}(#1)}
\newcommand{\C}[1]{{\rm C}(#1)}
\newcommand{\CFa}[1]{{\rm C_{F}^{ac}}(#1)}
\newcommand{\Cz}[1]{{\rm C}^+(#1)}
\newcommand{\Cf}[1]{{\rm C}^-(#1)}
\newcommand{\Cb}[1]{{\rm C^b}(#1)}
\newcommand{\K}[1]{{\rm K}(#1)}
\newcommand{\KF}[1]{{\rm K_{F}}(#1)}
\newcommand{\Ka}[1]{{\rm K^{ac}}(#1)}
\newcommand{\KFa}[1]{{\rm K_{F}^{ac}}(#1)}
\newcommand{\Kz}[1]{{\rm K}^+(#1)}
\newcommand{\Kbz}[1]{{\rm K}^{-,b}(#1)}
\newcommand{\Kbf}[1]{{\rm K}^{+,b}(#1)}
\newcommand{\KFb}[1]{{\rm K}^{-,Fb}(#1)}
\newcommand{\Kf}[1]{{\rm K}^-(#1)}
\newcommand{\Kb}[1]{{\rm K^b}(#1)}
\newcommand{\modcat}[1]{#1\mbox{{\rm -mod}}}
\newcommand{\smod}[1]{#1\mbox{{\rm -\underline{mod}}}}
\newcommand{\opp}{^{\rm op}}
\newcommand{\otimesL}{\otimes^{\rm\bf L}}
\newcommand{\rHom}{{\rm\bf R}{\rm Hom}\,}
\newcommand{\projdim}{\pd}
\newcommand{\lra}{\longrightarrow}
\newcommand{\ra}{\rightarrow}
\newcommand{\rad}{{\rm rad \, }}
\newcommand{\Lra}{\Longrightarrow}
\newcommand{\tra}{\twoheadrightarrow}
\newcommand{\Htp}{{\rm Htp}}
\newcommand{\Mod}{{\rm Mod}}
{\Large \bf
\begin{center}  Relative derived equivalences and relative
homological dimensions
 \end{center}}

\medskip

\centerline{\sc Shengyong Pan$^{1,2}$}

\begin{center} $^1$ Department of Mathematics,
 Beijing Jiaotong University,  Beijing 100044,
People's Republic of China
\end{center}

\begin{center} $^2$ School of Mathematical Sciences,
Beijing Normal University, Beijing 100875, People's Republic of
China\\ E-mail:shypan@bjtu.edu.cn
\end{center}

\medskip
\abstract{ Let $\mathscr{A}$ be a small abelian category. For a
closed subbifunctor $F$ of $\Ext_{\mathscr{A}}^{1}(-,-)$, Buan has
generalized the construction of the Verdier's quotient category to
get a relative derived category, where he localized with respect to
$F$-acyclic complexes. In this paper, the homological properties of
relative derived categories are discussed, and the relation with
derived categories is given. For Artin algebras, using relatively
derived categories, we give a relative version on derived
equivalences induced by $F$-tilting complexes. We discuss the
relationships between relative homological dimensions and relative
derived equivalences.}

\medskip {\small {\it 2000 AMS Classification}: 18E30,16G10;16S10,18G15.

\medskip {\it Key words:} relative derived category,
$F$-tilting complex, relative derived equivalence, relative
homological dimension.}

\section{Introduction}
Hochschild \cite{Hoch} introduced relative homological algebra in
categories of modules. And later, Heller, Butler and Horrocks
developed it in more general categories with a relative abelian
structure. Auslander and Solberg \cite{ASo1,ASo2,ASo3,ASo4,ASo5}
applied relative homological algebra to the representation theory of
Artin algebra. They studied relative homology in terms of
subbifunctors of the functor $\Ext^1(-,-)$ and developed the general
theory of relative cotilting modules for Artin algebras.

Derived categories were invented by Grothendieck-Verdier \cite{Ver}
in the early sixties. Today, they have widely been used in many
branches: algebraic geometry, stable homotopy theory, representation
theory, etc. In representation theory of Artin algebras, it is of
interest to investigate whether two Artin algebras have equivalent
derived categories. As is known, a Morita theory for derived
equivalences was established by Rickard \cite{Ri1}. According to his
theorem, two Artin algebras $A$ and $B$ are derived equivalent if
and only if there is a tilting complex for $A$ such that $B$ is
isomorphic to the endomorphism algebra of this complex.  An
interesting thing is to construct a new derived equivalence from a
given one by finding a tilting complex.
Rickard \cite{Ri2, Ri3} got
new derived equivalence by tensor products and trivial extensions.
In the recent years, Hu and Xi \cite{HX1,HX3} have provided various
techniques to construct new derived equivalences by finding tilting
complexes. Another interesting thing is to find some invariance
under derived equivalences. Derived equivalences between finite
dimensional algebras over a field have many invariants. For
instance, finiteness of finitistic dimension \cite{PX}, Hochschild
homology and cohomology \cite{Ri3} and $K$-theory \cite{DS} have
been shown to be invariant under derived equivalences.

Derived categories have been used effectively in relative
homological algebra. The main idea of relative homological algebra
is to replace projective modules by relative projective modules. It
is natural to study the corresponding version of the derived
category in this context. Since then, the relative derived
categories and relative tilting theory of Artin algebras have been
extensively studied. Recently, Gao and Zhang \cite{GZ} used
Gorenstein homological algebra to get Gorenstein derived categories.
Buan \cite{Bu1} also studied relative derived categories by
localizing relative quasi-isomorphisms. Using the notion of relative
derived categories, he generalized Happel's result on derived
equivalences induced by tilting module to the relative setting. Is
it very likely that relative tiltings provide something called
relative derived equivalences?

Motivated by Buan's work, we introduce relative derived equivalences
for Artin algebras in this paper. A 'Morita' theory for relative
derived categories is built, and some invariance of relative derived
equivalences is founded. The aim of this paper is to discuss the
relationships between relative derived equivalences and relative
homological dimensions.

To describe the main result, it is convenient to fix some notations.
Let $\Lambda$ and $\Gamma$ be Artin algebras. Assume that $F$ is a
subbifunctor of $\End^{1}_{\Lambda}(-,-)$ which is of finite type.
Denote by $\gldim(\Lambda)$ and $\fd(\Lambda)$ the global and the
finitistic dimensions of $\Lambda$, respectively. Denote by
$\gldim_{F}(\Lambda)$ and $\fd_{F}(\Lambda)$ the $F$-global and the
$F$-finitistic dimensions of $\Lambda$, respectively. Let $n\geq0$
be an integer. If the complex $\cpx{X}$ has the form:
$$
\cdots\ra X^{-n}\ra X^{-n+1}\ra\cdots\ra X^{-1}\ra X^{0}\ra\cdots,
$$
with $X^{i}\neq0$ and the differential being radical map for $-n\leq
i \leq0$, then $n$ is called the term length of $\cpx{X}$, denoted
by $t(\cpx{X})$.

Our main result can be stated as follows:

\noindent{\bf Theorem} (see Theorem 7.3) {\it Let
$L:\DFb{\Lambda}\ra \Db{\Gamma}$ be a relative derived equivalence.
Suppose $T^{\bullet}$ is an $F$-tilting complex for $\Lambda$
associated to $L$. Then

$(1)$ $\gldim_{F}(\Lambda)-t(T^{\bullet})\leq \gldim(\Gamma)\leq
\gldim_{F}(\Lambda)+t(T^{\bullet})+2$.

$(2)$ $\fd_{F}(\Lambda)-t(T^{\bullet})\leq \fd(\Gamma)\leq
\fd_{F}(\Lambda)+t(T^{\bullet})+2$. }

We give the upper and lower bounds of $\gldim(\Gamma)$ (resp.
$\fd(\Gamma)$) in the term length of relative tilting complex and
$\gldim_{F}(\Lambda)$ (resp. $\fd_{F}(\Lambda)$). In this theorem,
if $F=\Ext_{\Lambda}^{1}(-,-)$, then $T^{\bullet}$ is a tilting
complex for $\Lambda$ such that $\End(\cpx{T})\simeq\Gamma$. In the
proof of this theorem, we infer that $\gldim(\Lambda)$ (resp.
$\fd(\Lambda)$) and $\gldim(\Gamma)$ (resp. $\fd(\Gamma)$) satisfy a
similar formula but not the same as in this theorem, namely

$(1)$ $\gldim(\Lambda)-t(T^{\bullet})\leq \gldim(\Gamma)\leq
\gldim(\Lambda)+t(T^{\bullet})$.

$(2)$ $\fd(\Lambda)-t(T^{\bullet})\leq \fd(\Gamma)\leq
\fd(\Lambda)+t(T^{\bullet})$.

This paper is organized as follows. In Section 2, we consider closed
subbifunctors $F$ of $\Ext_{\mathscr{A}}^{1}(-,-)$, where
$\mathscr{A}$ is a small abelian category. We introduce the notion
of the relative derived category of $\mathscr{A}$ which was defined
in \cite{Bu1}. The homological properties of relative derived
categories are discussed, and the relation with derived categories
is given in Section 3. In Section 4, we give the triangulated
structure in the relative derived category. In Section 5, we show
that, the quotient of the bounded relative derived category by the
bounded homotopy category of relative projective complexes is
equivalent to the relative stable category as triangulated
categories. In Section 6, we introduce the notion of relative
derived equivalences for Artin algebras. In Section 7, we prove the
main result.

\section{Relative derived categories for abelian categories}
Let us explain the notion of relative derived categories. The notion
of relative derived categories was introduced earlier by Generalov
\cite {Ge}.

Let $\mathscr{A}$ be an abelian category. Suppose $A,
C\in\mathscr{A}$. Denote by $\Ext^{1}_{\mathscr{A}}(C,A)$ the set of
all exact sequences $0\ra A \ra B \ra C \ra 0 $ in $\mathscr{A}$
modulo the equivalence relation which is defined in the following
usual way. Two exact sequences are equivalent if the following
commutative diagram is commutative.
$$\xymatrix{
0\ar[r]&A\ar[r]\ar@{=}[d]   &    B\ar[r]\ar[d]    &    C\ar[r]\ar@{=}[d]& 0\\
0\ar[r]&A\ar[r]         &    B'\ar[r]          & C\ar[r]& 0}$$

Since that any additive category has finite direct sums and in
particular uniquely defined diagonal and codiagonal maps and that an
abelian category has pullback and pushout pairs, it follows that
$\Ext^{1}_{\mathscr{A}}(C,A)$ becomes an abelian group under Baer
sum. Therefore, $\Ext^{1}_{\mathscr{A}}(-,-)$ defines an additive
bifunctor $\mathscr{A}^{op} \times \mathscr{A} \lra \textbf{Ab}$,
where $\textbf{Ab}$ is the category of abelian groups.

Consider additive non-zero subbifunctors $F$ of $\Ext^{1}(-,-)$. To
each subfunctor corresponds a class of short exact sequences which
are called $F$-exact sequences. The class of $F$-exact sequences is
closed under the operations of pushout, pullback, Baer sums and
direct sums (see \cite{ASo1} or \cite{DRSS}). Given a subbifunctor
$F$ of $\Ext^{1}(-,-)$, we say that an exact sequence
$$\eta:0\ra X \ra Y \ra Z \ra 0 $$ in $\mathscr{A}$ is an $F$-exact
sequence if $\eta$ is in $F(Z,X)$. If $ 0\ra X \stackrel{f}\ra Y
\stackrel{g}\ra Z\ra 0 $ is an $F$-exact sequence, then $f$ is
called an $F$-monomorphism and $g$ is called an $F$-epimorphism. An
object $P$ is said to be $F$-projective if for each $F$-exact
sequence $0\ra X \ra Y \ra Z \ra 0 $, the sequence
$$0\ra \mathscr{A}(P,X) \ra \mathscr{A}(P,Y) \ra \mathscr{A}(P,Z)\ra
0$$ is exact. An object $I$ is called $F$-injective if for each
$F$-exact sequence $0\ra X \ra Y \ra Z \ra 0 $, the sequence $$0\ra
\mathscr{A}(X,I) \ra \mathscr{A}(Y,I) \ra \mathscr{A}(Z,I)\ra 0$$ is
exact. The subcategory of $\mathscr{A}$ consisting of all
$F$-projective (resp. $F$-injective) modules is denoted by
$\mathcal{P}(F)$ (resp. $\mathcal{I}(F)$).

We have the following characterization of when subbifunctors of
$\Ext^{1}(-,-)$ have enough projectives or injectives.

\begin{Lem} \label{2.1}$\rm\cite[Theorem\,1.12]{ASo1}$ Let $F$ be a subbifunctor of $\Ext^{1}(-,-)$.

$(1)$ $F$ has enough projectives if and only if $F=F_{\mathcal
{P}(F)}$ and $\mathcal {P}(F)$ is contravariantly finite in
$\Lambda$-mod.

$(2)$ $F$ has enough injectives if and only if $F=F^{\mathcal
{I}(F)}$ and $\mathcal {I}(F)$ is covariantly finite in
$\Lambda$-mod.

$(3)$ If there is a finite number of indecomposable relative
projectives (injectives) up to isomorphism, then there is also a
finite number of relative injectives (projectives), and these
numbers are the same.
\end{Lem}

Recall from \cite{Bu1, DRSS} that an additive subbifunctor $F$ is
said to be closed if the following equivalence statements hold.

1) The composition of $F$-epimorphisms is an $F$-epimorphism.

2) The composition of $F$-monomorphisms is an $F$-monomorphism.

3) For each object $X$ the functor $F(X,-)$ is half exact on
$F$-exact sequences.

4) For each object $X$ the functor $F(-,X)$ is half exact on
$F$-exact sequences.

5) The category $\mathscr{A}$ with respect to the $F$-exact
sequences is an exact category.

We will give some basic examples of closed subbifunctors. Let
$\mathcal {X}$ be a full subcategory of $\mathscr{A}$ and for each
pair of objects $A$ and $C$ in $\mathscr{A}$, we define
$$
F_{\mathcal {X}}(C,A)=\{0\ra A \ra B\ra C\ra 0\mid (\mathcal
{X},B)\ra(\mathcal {X},C)\ra 0 \text{\;is exact} \}.
$$
Dually we define for each pair of objects $A$ and $C$ in
$\mathscr{A}$
$$
F^{\mathcal {X}}(C,A)=\{0\ra A \ra B\ra C\ra 0\mid (B,\mathcal
{X})\ra(A,\mathcal {X})\ra 0 \text{\;is exact} \}.
$$

\begin{Lem}\label{2.2} $\rm\cite[Proposition\,1.7]{DRSS}$
The additive subbifunctors $F_{\mathcal {X}}$ and $F^{\mathcal {X}}$
of $\Ext^{1}_{\mathscr{A}}(-,-)$ are closed for any subcategory
$\mathcal {X}$ of $\mathscr{A}$.
\end{Lem}

\noindent{\bf Remark.} By Lemma \ref{2.1} and Lemma \ref{2.2}, it
follows that if a subbifunctor $F$ has enough projectives or
injectives, then $F$ is a closed subfunctor.

Let $\K{\mathscr{A}}$ be the homotopy category of complexes over
$\mathscr{A}$. We denote by $[1]$ the shift functor. A complex
$X^{\bullet}$ with differential $d_{X}^{\bullet}$ is said to be
$F$-acyclic if for each $i$ the induced complex $$0\ra Imd^{i-1} \ra
X^{i} \ra Imd^{i} \ra 0 $$ is an $F$-exact sequence. A map $h$ in
$\K{\mathscr{A}}$ is called an $F$-quasi-isomorphism if the mapping
cone $M(h)$ is an $F$-acyclic complex. If $F$ has enough
projectives, then $F$-exact sequences are exact sequences and
$F$-acyclic complexes are acyclic complexes. But the converse is not
true.

If the class of $F$-acyclic complexes are closed under the operation
of mapping cone, then we call $F$ a triangulated subbifunctor. In
this case the class of $F$-acyclic complexes is a null-system in
$\K{\mathscr{A}}$. Since the class of $F$-quasi-isomorphisms is a
multiplicative system, it follows that we localize with respect to
this system.

\begin{Lem} $\rm\cite[Theorem\,2.4]{Bu1}$
A subbifunctor is triangulated if and only if it is closed.
\end{Lem}

Now we assume that $F$ is a closed subbifunctor. It is clear that
$$\mathscr{N}=\{ X^{\bullet}\in \text{Obj}(\K{\mathscr{A}})\mid
X^{\bullet} \text{\;is an F-acyclic complex} \}$$ is a null system
and a thick subcategory of $\K{\mathscr{A}}$. Define a morphism set
\begin{eqnarray*}
\Sigma(\mathscr{N})=\{X^{\bullet} \stackrel{f^{\bullet}}\ra
Y^{\bullet} \mid\text{such that\;} X^{\bullet}
\stackrel{f^{\bullet}}\ra Y^{\bullet} \ra Z^{\bullet} \ra
X^{\bullet}[1]\text{\; is a distinguished triangle
in\;}\\\K{\mathscr{A}} \text{\;with\;} Z^{\bullet} \in\mathscr{N}
\}.
\end{eqnarray*} The relative derived
category of $\mathscr{A}$ is defined to be the Verdier quotient,
that is,
$\DFb{\mathscr{A}}:=\K{\mathscr{A}}/\mathscr{N}=\Sigma(\mathscr{N})^{-1}\K{\mathscr{A}}$.

The objects of $\DFb{\Lambda}$ are the same as for
$\K{\mathscr{A}}$. A map in $\DFb{\Lambda}$ from $X^{\bullet}$ to
$Y^{\bullet}$ is the equivalence class of "roofs", that is, of
fractions $\cpx{g}/\cpx{f}$ of the form $\xymatrix{ \cpx{X}&
Z^{\bullet} \ar_{f^{\bullet}}[l]\ar^{g^{\bullet}}[r]& \cpx{Y}}$,
where $Z^{\bullet}\in \K{\mathscr{A}}$, $f^{\bullet}:Z^{\bullet}\ra
X^{\bullet} $ is an $F$-quasi-isomorphism, and
$g^{\bullet}:Z^{\bullet}\ra X^{\bullet}$ is a morphism in
$\K{\mathscr{A}}$. Two such roofs $\xymatrix{ \cpx{X}&  Z^{\bullet}
\ar_{f^{\bullet}}[l]\ar^{g^{\bullet}}[r]& \cpx{Y}}$ and $\xymatrix{
\cpx{X'}&  \cpx{Z'}\ar_{\cpx{f'}}[l]\ar^{\cpx{g'}}[r]& \cpx{Y'}}$
are equivalent if there exists a commutative diagram
$$\xymatrix{
   &  Z^{\bullet} \ar_{f^{\bullet}}[dl]\ar^{g^{\bullet}}[dr]        \\
X^{\bullet}   & W^{\bullet}\ar_{h^{\bullet}}[l]\ar[u]\ar[d]&    Y^{\bullet} , \\
&  Z'^{\bullet} \ar^{f'^{\bullet}}[ul]\ar_{g'^{\bullet}}[ur] }$$
where $h^{\bullet}$ is an $F$-quasi-isomorphism. Note that the
diagram of the form $\xymatrix{ \cpx{X}&  Z^{\bullet}
\ar_{f^{\bullet}}[l]\ar^{g^{\bullet}}[r]& \cpx{Y}}$ will be called a
left roof.

\section{Homological properties of relative
derived categories}

In this section, we give some basic properties of relative derived
categories.

The following lemma is well-known.

\begin{Lem}$\rm\cite[Lemma\,2.4]{GZ}$
Let $\mathcal {T}_{1}$ and $\mathcal {T}_{2}$ be triangulated
subcategories of a triangulated category $\mathcal {T}$, and
$\mathcal {T}_{2}$ a full triangulated subcategory of $\mathcal
{T}_{1}$. Then there is an isomorphism of triangulated categories
$(\mathcal {T}/\mathcal {T}_{2})/(\mathcal {T}_{1}/\mathcal {T}_{2})
\simeq \mathcal {T}/\mathcal {T}_{1}$.
\end{Lem}

Set $\mathscr{N'}$=$\{ X^{\bullet}\mid X^{\bullet}$ is an acyclic
complex $\}$. Then $\mathscr{N'}$ is also a null system.
Consequently, $\mathscr {N}$ is a triangulated subcategory of
$\mathscr{N'}$. By Lemma 3.1, it follows that
$$
\Db{\mathscr{A}}=\K{\mathscr{A}}/\mathscr{N'}\simeq
(\K{\mathscr{A}}/\mathscr{N})/(\mathscr{N'}/\mathscr{N})\simeq
\DFb{\mathscr{A}}/(\mathscr{N'}/\mathscr{N}),
$$
where $\Db{\mathscr{A}}$ is the derived category of $\mathscr{A}$.
Then the usual derived category is a quotient category of relative
derived category.

As is known, if $F$ has enough projective objects, then an exact
sequence $0\ra A\ra B\ra C\ra 0$ is $F$-exact if and only if $0\ra
\Hom(P,A)\ra \Hom(P,B)\ra \Hom(P,C)\ra 0$ is exact for all $P$ in
$\mathcal {P}(F)$. We have the following lemma.

\begin{Lem} If $F$ has enough projective objects, then a complex $X^{\bullet}$
is an $F$-acyclic complex if and only if $\Hom(P,X^{\bullet})$ is an
acyclic complex for all $P$ in $\mathcal {P}(F)$.
\end{Lem}

It is well known that the categories $\D^{-}(\mathscr{A})$,
$\D^{+}(\mathscr{A})$ and $\Db{\mathscr{A}}$ are full subcategories
of $\D(\mathscr{A})$. This remains true for
$\D^{\ast}_{F}(\mathscr{A})$ in $\D_{F}(\mathscr{A})$, where
$\star=+,-,b$. We will prove these from a series of observations to
realize these facts.

\begin{Prop}$\rm\cite[Lemma\,10.1.13]{W}$ Let $\mathcal {S}$ be
a multiplicative system of $\mathcal {C}$ and $\mathcal {D}$ a full
subcategory of $\mathcal {C}$. Assume that $\mathcal {S}\cap\mathcal
{D}$ is a multiplicative system of $\mathcal {D}$, and one of the
following conditions is satisfied:

$(1)$ If $s: X'\ra X$ is a morphism in $\mathcal {S}$, with
$X\in\mathcal {D}$, then there is a morphism $f:X^{''}\ra X'$ in
$\mathcal {C}$ such that $X^{''}\in \mathcal {D}$ and $fs\in\mathcal
{S}$.

$(2)$ If $s: X\ra X'$ is a morphism in $\mathcal {S}$ with
$X\in\mathcal {D}$, then there is a morphism $f: X'\ra X^{''}$ in
$\mathcal {C}$ such that $X^{''}\in \mathcal {D}$ and $sf\in\mathcal
{S}$.

Then the natural functor $(\mathcal {S}\cap\mathcal
{D})^{-1}\mathcal {D}\ra\mathcal {S}^{-1}\mathcal {C}$ is fully
faithful. That is, $(\mathcal {S}\cap\mathcal {D})^{-1}\mathcal {D}$
is considered as a full subcategory of $\mathcal {S}^{-1}\mathcal
{C}$.
\end{Prop}

\begin{Koro}
$(1)$ If $X'^{\bullet}\ra X^{\bullet}$ is an $F$-quasi-isomorphism,
$X^{\bullet}\in \Kb{\mathscr{A}}$, $X'^{\bullet}\in
\Kz{\mathscr{A}}$, then there is an $F$-quasi-isomorphism
$X^{''}{^{\bullet}}\ra X'^{\bullet}$ such that
$X^{''}{^{\bullet}}\in \Kb{\mathscr{A}}$.

$(2)$ If $Y^{\bullet}\ra X^{\bullet}$ is an $F$-quasi-isomorphism,
$X^{\bullet}\in \Kf{\mathscr{A}}$, $Y^{\bullet}\in \K{\mathscr{A}}$,
then there is an $F$-quasi-isomorphism $Y'^{\bullet}\ra Y^{\bullet}$
such that $Y'^{\bullet}\in \Kf{\mathscr{A}}$.
\end{Koro}

\textbf{Proof.} The proof of ($1$) is similar to the one of ($2$).
We only need to prove the second case.

(2) Let $X^{\bullet}$ be in $\Kf{\mathscr{A}}$, we may assume that
$X^{i}=0$ for $i>0$. Then $X^{\bullet}$ has the form
$$
\cdots\ra X^{-1}\ra X^{0}\ra0\ra \cdots.
$$
Let $f^{\bullet}: X^{\bullet}\ra Y^{\bullet}$ be a map of complexes.
Then we have the following commutative diagram
$$\xymatrix{
\cdots \ar[r]&Y^{-1}\ar^{d_{Y}^{-1}}[r]\ar^{f^{-1}}[d]&Y^{0}
\ar^{d_{Y}^{0}}[r]\ar^{f^{0}}[d] & Y^{1}\ar[r]\ar[d]& \cdots\\
\cdots \ar[r]&X^{-1}\ar^{d^{-1}}[r] &  X^{0}\ar[r]  & 0\ar[r]&
\cdots.}
$$
Since $f^{\bullet}$ is an $F$-quasi-isomorphism, it follows that
$M(f^{\bullet})$ is an $F$-acyclic complex.
$$
M(f^{\bullet}):\cdots\ra X^{-1}\oplus Y^{0}\ra X^{0}\oplus Y^{1}\ra
Y^{2}\ra \cdots
$$
Decompose $M(f^{\bullet})$ into the direct sum of two $F$-exact
complexes:
$$
\cdots\ra X^{-1}\oplus Y^{0}\ra X^{0}\oplus Y^{1}\ra
Imd^{1}_{Y}\ra0\ra \cdots
$$
and
$$
\cdots\ra 0\ra Imd_{Y}^{2}\ra  Y^{3}\ra  \cdots.
$$
Then we get the following commutative diagram:
$$\xymatrix{
\cdots \ar[r]&Y^{0}\ar[r]\ar@{=}[d]   &    Y^{1}\ar[r]\ar@{=}[d]
&  Imd^{1}_{Y}\ar[r]\ar[d]& 0 \ar[r]\ar[d] &\cdots\\
\cdots \ar[r]&Y^{0}\ar[r]\ar[d]   & Y^{1}\ar[r]\ar[d] &
Y^{2}\ar[r]\ar[d]& Y^{3}\ar[r]\ar@{=}[d]& \cdots \\
\cdots \ar[r]&0\ar[r] & 0\ar[r] & Imd_{Y}^{2}\ar[r]& Y^{3}\ar[r]&
\cdots.}$$ Set
$$
Y'^{\bullet}:\cdots\ra Y^{0}\ra Y^{1}\ra Imd^{1}\ra 0\ra \cdots.
$$
Then $Y'^{\bullet}$ and $Y^{\bullet}$ are $F$-quasi-isomorphic with
$Y'^{\bullet}\in \Kf{\mathscr{A}}$.$\square$

\begin{Prop}
$\DFb{\mathscr{A}}$, $\DFf{\mathscr{A}}$ and $\DFz{\mathscr{A}}$ are
full subcategories of $\DF{\mathscr{A}}$.
\end{Prop}

\textbf{Proof.} By definition,
$\DFf{\mathscr{A}}=\Kf{\mathscr{A}}/(\mathscr{N}\cap
\Kf{\mathscr{A}})$, since $\Kf{\mathscr{A}}$ is a full subcategory
of $\K{\mathscr{A}}$. Moreover, if $Y^{\bullet}\ra X^{\bullet}$ is
an $F$-quasi-isomorphism, $X^{\bullet}\in \Kf{\mathscr{A}}$,
$Y^{\bullet}\in \K{\mathscr{A}}$, then there is an
$F$-quasi-isomorphism $Y'^{\bullet}\ra X^{\bullet}$ by Corollary
3.4(2). From Proposition 3.3, we deduce that $\DFf{\mathscr{A}}$ is
a full subcategory of $\DF{\mathscr{A}}$. In a similar way,
$\DFf{\mathscr{A}}$ and $\DFz{\mathscr{A}}$ are full subcategories
of $\DF{\mathscr{A}}$.

The following results make the morphisms in $\DF{\mathscr{A}}$ to
the morphisms in $\K{\mathscr{A}}$ to understand.

\begin{Lem}Let $P^{\bullet}\in \Kf{\mathcal {P}(F)}$ and
$f^{\bullet}:X^{\bullet}\ra P^{\bullet}$ an $F$-quasi-isomorphism.
Then there is a map $g^{\bullet}:P^{\bullet}\ra X^{\bullet}$ such
that $g^{\bullet} f^{\bullet}$ is homotopic to $id_{P^{\bullet}}$.
\end{Lem}

\textbf{Proof.} Let $f^{\bullet}:X^{\bullet}\ra P^{\bullet}$ be an
$F$-quasi-isomorphism. Then we have a distinguished triangle
$X^{\bullet}\stackrel{f^{\bullet}}\ra
P^{\bullet}\stackrel{h^{\bullet}}\ra
M(f^{\bullet})\stackrel{l^{\bullet}}\ra X^{\bullet}[1]$ in
$\Kf{\mathcal {P}(F)}$, where the mapping cone
$M(f^{\bullet}):=X^{\bullet}[1]\oplus P^{\bullet}$ is an $F$-acyclic
complex. For $P^{\bullet}\in \Kf{\mathcal {P}(F)}$, without loss of
generality we can assume that $P^{i}=0$ for $i>n$. Therefore, we
have the following commutative diagram
$$\xymatrix{
\cdots \ar[r]&X^{n-2}\ar^{d_{X}^{n-2}}[r]\ar_{f^{n-2}}[d] &
X^{n-1}\ar^{d_{X}^{n-1}}[r]\ar_{f^{n-1}}[d] &
 X^{n}\ar^{d_{X}^{n}}[r]\ar_{f^{n}}[d]& X^{n+1} \ar[r]\ar[d] &\cdots\\
\cdots \ar[r]&P^{n-2}\ar[r]\ar_{h^{n-2}}[d] &
P^{n-1}\ar[r]\ar@{-->}_{s^{n-1}}[dl]\ar_{h^{n-1}}[d] &
P^{n}\ar[r]\ar@{-->}_{s^{n}}[dl]\ar_{h^{n}}[d]& 0\ar[r]\ar[d]& \cdots \\
\cdots \ar[r]&M(f^{\bullet})^{n-2}\ar[r]&
M(f^{\bullet})^{n-1}\ar^{d^{n-1}}[r] &
M(f^{\bullet})^{n}\ar^{d^{n}}[r]& M(f^{\bullet})^{n+1}\ar[r]&
\cdots.}$$ From the above diagram, we have
$h^{n}d_{M(f^{\bullet})}^{n}=0$. By the (relative) Comparision
Theorem, it follows that $h^{\bullet}$ is a null-homotopic. Then
there is a chain map
$s^{\bullet}=(g^{\bullet},v^{\bullet}):P^{\bullet}\ra
X^{\bullet}\oplus P^{\bullet}[1]$ such that $h^{\bullet}=s^{\bullet}
d_{M(f^{\bullet})}^{\bullet}+d_{P}^{\bullet} s^{\bullet}$, where
$d_{M(f^{\bullet})}^{i}=\left(\begin{array}{cc}
-d_X^{i+1}&f^{i+1}\\
0&d_P^i
\end{array}\right)$. Thus the above equation yields the following
equations
$$
-gd_X+d_Pg=0, \; gf+vd_P+d_Pv=id.
$$
It follows that $g$ is a chain map and $id=g^{\bullet} f^{\bullet}$
in $\Kf{\mathscr{A}}$.

\begin{Prop} \label{3.7}Let $P^{\bullet}\in \Kf{\mathcal {P}(F)}$ and $X^{\bullet}\in \K{\mathscr{A}}$.
Then
$$
\Hom_{\K{\mathscr{A}}}(P^{\bullet},X^{\bullet})\simeq
\Hom_{\DF{\mathscr{A}}}(P^{\bullet},X^{\bullet}).
$$
\end{Prop}

\textbf{Proof.} Define a map $\varphi:
\Hom_{\K{\mathscr{A}}}(P^{\bullet},X^{\bullet})\ra
\Hom_{\DF{\mathscr{A}}}(P^{\bullet},X^{\bullet})$ which sends
$f^{\bullet}$ to $f^{\bullet}/id$. If $f^{\bullet}: P^{\bullet}\ra
X^{\bullet}$ satisfies $f^{\bullet}/id = 0$, then there is an
$F$-quasi-isomorphism $t^{\bullet}: X'^{\bullet}\ra P^{\bullet}$
such that $t^{\bullet} f^{\bullet}$ is homotopic to zero. By Lemma
3.6, there is an $F$-quasi-isomorphism $g^{\bullet}:P^{\bullet}\ra
X'^{\bullet}$ such that $g^{\bullet} t^{\bullet}$ is homotopic to
$id_{P^{\bullet}}$. Thus, we get $g^{\bullet} t^{\bullet}
f^{\bullet}=0$. Consequently, $f^{\bullet}=0$ in $\K{\mathscr{A}}$.
This implies that $\varphi$ is an injective map. For any
$f^{\bullet}/s^{\bullet}$ is in
$\Hom_{\DF{\mathscr{A}}}(P^{\bullet},X^{\bullet})$, by Lemma 3.6,
there is an $F$-quasi-isomorphism $g^{\bullet}$ such that
$g^{\bullet} s^{\bullet}$ is homotopic to $id_{P^{\bullet}}$.
Consequently, $f^{\bullet}/s^{\bullet} =
f^{\bullet}g^{\bullet}/id_{P^{\bullet}}$. Therefore, $\varphi$ is
surjective. $\square$

Suppose that the functor $F$ has enough projective objects. For any
$X$ in $\mathscr{A}$, take a minimal $F$-projective resolution
$P_{X}$ of $X$, that is, there exists an $F$-exact sequence
$$
\cdots\ra P^{-2}\ra P^{-1}\ra P^{0}\ra X\ra 0,
$$with $P^{i}$ $F$-projectives for $i>0$.
Then we define $\Ext^{i}_{F}(X,Y)$ to be the $i$-th homology of the
complex
$$
0\ra \Hom(P^{0},Y)\ra \Hom(P^{-1},Y)\ra\cdots \ra\Hom(P^{-i},Y)\ra
\cdots.
$$

\begin{Lem} Suppose that $F$ has enough projective objects. Let $X$ and $Y$
be in $\mathscr{A}$. Then
$$
\Ext^{i}_{F}(X,Y)\simeq\Hom_{\DFb{\mathscr{A}}}(X,Y[i])
$$
for all $i\in \mathbb{Z}$.
\end{Lem}

\textbf{Proof.} Take a minimal $F$-projective resolution
$P_{X}^{\bullet}$ of $X$ with $P_{X}^{\bullet}\in
\Kbz{\mathscr{A}}$. Then we have the following
\begin{eqnarray*}
\Hom_{\DFb{\mathscr{A}}}(X,Y[i])\simeq
\Hom_{\DFb{\mathscr{A}}}(P_{X}^{\bullet},Y[i])\simeq
\Hom_{\K{\mathscr{A}}}(P_{X}^{\bullet},Y[i])\\\simeq
H^{i}(\Hom(P_{X},Y)) = \Ext_{F}^{i}(X,Y).
\end{eqnarray*}

Define $\KFb{\mathcal {P}(F)}$ to be the subcategory of
$\Kbz{\mathcal {P}(F)}$ by
$$
\KFb{\mathcal {P}(F)}=\{\cpx{X}\in \Kf{\mathcal {P}(F)}| \exists
\;n\in\mathbb{Z}\text{\; such that \;} H^i(\Hom_A(P,\cpx{X}))=0
\text{\; for\;} i\leq n, \forall P\in \mathcal {P}(F)\}.
$$
It is easy to easy that for $\cpx{X}\in \Kf{\mathcal {P}(F)}$ and
$\forall i\in \mathbb{Z}$, if $H^i(Hom_A(P,\cpx{X}))=0$, for $i\leq
n$, $\forall P\in \mathcal {P}(F)$, then $H^i(\cpx{X})=0$.

\begin{Prop}Suppose $F$ has enough projective objects. Then the
 natural functor induces two triangle equivalences:
 $$
 \DFf{\mathscr{A}}\simeq \Kf{\mathcal {P}(F)},{\;}
 \DFb{\mathscr{A}}\simeq \KFb{\mathcal {P}(F)}.
$$
\end{Prop}

\textbf{Proof.} For more details, we refer the reader to
\cite[Theorem 10.4.8]{W} and \cite[Proposition 6.3.1]{KZ}. $\square$

Let $H: \mathscr{A} \ra \DFb{\mathscr{A}}$ be the composition of the
embedding $\mathscr{A} \ra \Kb{\mathscr{A}}$ with the localization
functor $\Kb{\mathscr{A}}\ra \DFb{\mathscr{A}}$. Then we have the
following proposition.
\begin{Lem}
The functor $H: \mathscr{A} \ra \DFb{\mathscr{A}}$ is fully
faithful.
\end{Lem}
\textbf{Proof.} It suffices to show that $ H(X,Y):
\Hom_{\mathscr{A}}(X,Y)\simeq\Hom_{\DFb{\mathscr{A}}}(X,Y)$ as
abelian groups. Let $u\in\Hom_{\mathscr{A}}(X,Y)$ such that
$H(u)=0$. Then we have the following diagram
$$\xymatrix{
   &  X \ar@{=}[dl]\ar^{u}[dr]  \\
X   & M^{\bullet}\ar_{a^{\bullet}}[l]\ar[r]\ar_{a^{\bullet}}[u]\ar[d]&    Y , \\
&  X \ar@{=}[ul]\ar_{0}[ur] }$$ where $a^{\bullet}: M^{\bullet}\ra
X$ is an $F$-quasi-isomorphism such that $a^{\bullet} u=0$.
Therefore, $H^{0}(a):H^{0}(M^{\bullet})\ra X$ is an isomorphism and
$H^{0}(a) u=0$. This implies that $u=0$ and $H$ is a faithful
functor. We take a left roof $ \xymatrix{ X    &  U^{\bullet}
\ar_{s^{\bullet}}[l]\ar^{f^{\bullet}}[r] & Y } $ in
$\Hom_{\DFb{\mathscr{A}}}(X,Y)$, where $U^{\bullet}\in
\Kb{\mathscr{A}}$ and $s^{\bullet}$ is an $F$-quasi-isomorphism.
Consequently, $H^{0}(s^{\bullet}):H^{0}(U^{\bullet})\ra X$ is an
isomorphism in $\mathscr{A}$. Set $g=H^{0}(s^{\bullet})^{-1}
 H^{0}(f^{\bullet}): X\ra Y$. Consider the truncation complex
$W^{\bullet}:\cdots\ra U^{-2}\ra U^{-1}\ra Kerd^{0}\ra0$ of
$U^{\bullet}$ and the canonical chain map
$i^{\bullet}:W^{\bullet}\ra U^{\bullet}$. It is easy to see that
$i^{\bullet}$ is an $F$-quasi-isomorphism. Since $s^{\bullet}$ is an
$F$-quasi-isomorphism, it follows that $i^{\bullet} s^{\bullet}$ is
an $F$-quasi-isomorphism. To get the following commutative diagram
$$\xymatrix{
   &  U^{\bullet} \ar_{s^{\bullet}}[dl]\ar^{f^{\bullet}}[dr]        \\
X   & W^{\bullet}\ar_{i^{\bullet}
s^{\bullet}}[l]\ar[r]\ar_{i^{\bullet}}[u]\ar[d]&    Y ,\\
&  X \ar@{=}[ul]\ar_{g}[ur] .}$$ it is sufficient to show that
$i^{\bullet} f^{\bullet}=i^{\bullet} s^{\bullet} g$. From the
following commutative diagram
$$\xymatrix{
 W^{\bullet}\ar[d] \ar[r]^{i^{\bullet}}
 &U^{\bullet} \ar[d]^{s^{\bullet}}  \\
  H^{0} (U^{\bullet})  \ar[r]^{H^{0}(s^{\bullet})}
  &   X ,} $$
we conclude that $i^{\bullet} s^{\bullet} g=i^{\bullet} s^{\bullet}
H^{0}(s^{\bullet})^{-1} H^{0}(f^{\bullet})=i^{\bullet} f^{\bullet}$.
Then $L(g/id)=f^{\bullet}/s^{\bullet}$ and $L$ is a full functor.
$\square$

\section{Triangles in relative derived categories}

In this section, we describe the triangulated structure of
$\DFb{\mathscr{A}}$, where $F$ has enough projectives.

We have the following proposition.

\begin{Prop}Suppose that
$$
0\ra X^{\bullet}\stackrel{f^{\bullet}}\ra Y^{\bullet}
\stackrel{g^{\bullet}}\ra Z^{\bullet} \ra 0
$$
is an $F$-exact sequence of complexes in $\C{\mathscr{A}}$, that is,
$$
0\ra X^{i} \ra Y^{i} \ra Z^{i} \ra 0
$$
is an $F$-exact sequence for each $i$. Let $M(f^{\bullet})$ be the
mapping cone of $f^{\bullet}$, and let $\phi^{i}:
M(f^{\bullet})^{i}=X^{i+1}\oplus Y^{i} \ra Z^{i}$ be the morphism
$\left(\begin{array}{cc}
0\\
g^i
\end{array}\right)$.
Then $\phi^{\bullet}: M(f^{\bullet})\ra Z^{\bullet}$ is a morphism
of complexes, $\alpha(f^{\bullet}) \phi^{\bullet}=g^{\bullet}$, and
$\cpx{\phi}$ is an $F$-quasi-isomorphism. Moreover,
$$
X^{\bullet} \ra Y^{\bullet} \ra Z^{\bullet} \ra X^{\bullet}[1]
$$
is a distinguished triangle in $\DFb{\mathscr{A}}$.
\end{Prop}

{\textbf{Proof}.} It is straightforward to see that $\phi^{\bullet}$
is a morphism of complexes. We know that there is a diagram of maps
of complexes:

$$\xymatrix{
X^{\bullet}\ar[r]&Y^{\bullet}\ar^{\alpha(f^{\bullet})}[r]\ar@{=}[d]
&
M(f^{\bullet})\ar^{\beta(f^{\bullet})}[r]\ar^{\phi^{\bullet}}[d]  & X^{\bullet}[1]\\
X^{\bullet}\ar[r]&Y^{\bullet}\ar[r]^{g^{\bullet}} &
Z^{\bullet}\ar[r] & X^{\bullet}[1].}$$ It suffices to show that
$M(f^{\bullet})\stackrel{\phi^{\bullet}}\ra Z^{\bullet}$ is an
$F$-quasi-isomorphism. Moreover, there is an $F$-exact sequence of
complexes $0\ra M(id_{\cpx{X}})\stackrel{\rho^{\bullet}}\ra
M(f^{\bullet})\stackrel{\phi^{\bullet}}\ra Z^{\bullet}\ra 0$, where
$\rho^{i}= \left(\begin{array}{cc}
                 id & 0\\ 0& f^i
                 \end{array}\right)$.
By assumption, for each $P\in\mathcal {P}(F)$, it follows that
 $$0\ra (P,M(id_{\cpx{X}})) \ra (P,M(f^{\bullet})) \ra (P,Z^{\bullet}) \ra 0$$
 is an exact sequence of complexes.
Since $H^{i}((P,M(id_{\cpx{X}})))=0$, we deduce that
$\phi^{\bullet}$ is an $F$-quasi-isomorphism. Therefore,
$\cpx{\phi}$ is an isomorphism in $\DFb{\mathscr{A}}$. Therefore,
$$
X^{\bullet} \stackrel{f^{\bullet}}\ra Y^{\bullet}
\stackrel{g^{\bullet}}\ra Z^{\bullet} \stackrel{h^{\bullet}} \ra
X^{\bullet}[1]
$$
is a distinguished triangle in $\DFb{\mathscr{A}}$, where
$h^{\bullet}=\phi^\bullet{}^{-1} \beta(f^{\bullet})$. $\square$

\section{Relative stable categories and relative derived categories }

In this section, we assume that $F$ has enough projectives and
injectives and $\mathcal {P}(F)=\mathcal {I}(F)$. Then $\mathscr{A}$
is a Frobenius category. Let $\mathscr{A}/\mathcal {P}(F)$ denote
the stable category of $\mathscr{A}$. By \cite[Theorem 2.6]{Ha1},
$\mathscr{A}/\mathcal {P}(F)$ is a triangulated category. The shift
functor of $\mathscr{A}/\mathcal {P}(F)$ is $\Omega^{-1}_{F}$, where
$\Omega_{F}$ is the $F$-syzygy functor. Since $\mathcal {P}(F)$ is
closed under taking direct summands, it follows that $\Kb{\mathcal
{P}(F)}$ is a thick subcategory of $\DFb{\mathscr{A}}$.

We are now ready to state our main result of this section.
\begin{Theo} The quotient category
$$
\DFb{\mathscr{A}}/\Kb{\mathcal {P}(F)}
$$ is equivalent as a triangulated
category to $\mathscr{A}/\mathcal {P}(F)$.
\end{Theo}

\textbf{Proof.} The proof is similar to that of  \cite[Theorem
2.1]{Ri2}. For convenience, we give more details.

Consider the additive functor $$L': \mathscr{A} \ra
\DFb{\mathscr{A}}/\Kb{\mathcal {P}(F)}$$ obtained by the composition
of the embedding $\mathscr{A} \ra \DFb{\mathscr{A}}$ with the
localization functor $\DFb{\mathscr{A}}\ra
\DFb{\mathscr{A}}/\Kb{\mathcal {P}(F)}$. The natural inclusion of
$\mathscr{A}$ sends $F$-projective objects to zero in the quotient
category, hence the functor $L$ factors through the stable category,
inducing an additive functor $L: \mathscr{A}/\mathcal {P}(F)\simeq
\DFb{\mathscr{A}}/\Kb{\mathcal {P}(F)}$.

The proof will be divided into three steps.

Step1. $L$ is an exact functor. Consider a distinguished triangle
$$X\ra Y\ra Z\ra \Omega^{-1}_{F}(X)$$ in $\mathscr{A}/\mathcal
{P}(F)$. This distinguished triangle is from the following pushout
diagram:
$$\xymatrix{
0 \ar[r]&X\ar[r]\ar[d]   &  I(X)  \ar[r]\ar[d]    &    \Omega^{-1}_{F}(X)\ar[r]\ar@{=}[d]& 0\\
0 \ar[r]&Y\ar[r]         &    Z\ar[r]          &
\Omega^{-1}_{F}(X)\ar[r]& 0.}$$ Since $F$-exact sequences are closed
under pushout, it follows that $0\ra Y\ra Z\ra
\Omega^{-1}_{F}(X)\ra0$ is an $F$-exact sequence. From $F$-exact
sequences  $0\ra X\ra I(X)\ra \Omega^{-1}_{F}(X)\ra 0$ and $0\ra
I(X)\ra Z\ra \Omega^{-1}_{F}(X)\ra 0$, we have distinguished
triangles $X\ra I(X)\ra \Omega^{-1}_{F}(X)\ra X[1]$ and $Y\ra Z\ra
\Omega^{-1}_{F}(X)\ra Y[1]$ in $\DFb{\mathscr{A}}$.

Since $L'(I(X))=0$ in $\DFb{\mathscr{A}}/\Kb{\mathcal {P}(F)}$, we
deduce that $$L(X)\ra L(Y)\ra L(Z)\ra L(X)[1]$$ is a distinguished
triangle in $\DFb{\mathscr{A}}/\Kb{\mathcal {P}(F)}$. Consequently,
$L$ is an exact functor.

Step2. $L$ is a fully faithful functor. Firstly, $L'$ is a full
functor. We see that the map
$$
\Hom_{\mathscr{A}}(X,Y)\ra \Hom_{\DFb{\mathscr{A}}/\Kb{\mathcal
{P}(F)}}(X,Y)
$$
sending $f: X\ra Y$ to $f/id$ is surjective, where $f/id$ is the
left roof
$$\xymatrix{
    X \ar_{id}@{=}[r]& X \ar^{f}[r] & Y  .}$$Take a left roof
$\xymatrix{
 X  &  \cpx{Z} \ar_{s}[l]\ar^{g}[r] & Y }$ in
$\Hom_{\DFb{\mathscr{A}}/\Kb{\mathcal {P}(F)}}(X,Y)$, where
$\cpx{Z}\in \DFb{\mathscr{A}}$. It follows that there is a
distinguished triangle
$$
\cpx{Z}\stackrel{s}\ra X\ra M(s)\ra \cpx{Z}[1]
$$
in $\DFb{\mathscr{A}}$ such that the mapping cone $M(s)$ is in
$\Kb{\mathcal {P}(F)}$. Applying the cohomological functor
$\Hom_{\DFb{\mathscr{A}}/\Kb{\mathcal {P}(F)}}(-,Y)$ to the
distinguished triangle $\cpx{Z}\stackrel{s}\ra X\ra M(s)\ra
\cpx{Z}[1]$, we have an exact sequence
\begin{eqnarray*}
\cdots\ra\Hom_{\DFb{\mathscr{A}}/\Kb{\mathcal {P}(F)}}(M(s),Y)\ra
\Hom_{\DFb{\mathscr{A}}/\Kb{\mathcal {P}(F)}}(X,Y)
\\\stackrel{(s,Y)}\ra \Hom_{\DFb{\mathscr{A}}/\Kb{\mathcal
{P}(F)}}(\cpx{Z},Y)\ra \Hom_{\DFb{\mathscr{A}}/\Kb{\mathcal
{P}(F)}}(M(s)[-1],Y)\ra \cdots
\end{eqnarray*}
Furthermore,
$$
\Hom_{\DFb{\mathscr{A}}/\Kb{\mathcal
{P}(F)}}(M(s)[-1],Y)=0\quad\text{since}\quad M(s)[-1]\quad
\text{is\; in}\quad \Kb{\mathcal {P}(F)}.$$ Then we get $g=s f$ and
the following commutative diagram
$$\xymatrix{
   &  \cpx{Z} \ar_{s}[dl]\ar^{g}[dr]        \\
X   & \cpx{Z} \ar^{s}[l]\ar[r]\ar@{=}[u]\ar^{s}[d]& Y .\\
&  X \ar@{=}[ul]\ar_{f}[ur] }$$ Therefore, $f/id$ and $g/s$ are
equivalent in $\DFb{\mathscr{A}}/\Kb{\mathcal {P}(F)}$. This implies
that $L$ is a fully functor.

We claim that if $X$ is an object of $\mathscr{A}$, then $L(X)=0$ if
and only if $X=0$. Indeed, no non-projective object is isomorphic to
an object of $\Kb{\mathcal {P}(F)}$.

Assume that $\alpha: X\ra Y$ is in $\mathscr{A}/\mathcal{P}(F)$ such
that $L(\alpha)=0$. There is a distinguished triangle
$X\stackrel{\alpha}\ra Y\stackrel{x}\ra Z\ra \Omega^{-1}_{F}(X)$ in
$\mathscr{A}/\mathcal {P}(F)$. Since $L$ is an exact functor, we see
that $L(X)\stackrel{L(\alpha)}\ra L(Y)\ra L(Z)\ra L(X)[1]$ is a
distinguished triangle in $\DFb{\mathscr{A}}/\Kb{\mathcal{P}(F)}$.
Since $L(\alpha)=0$, we conclude that $L(Z)\simeq (L(Y) \oplus
L(X)[1])$. Thus, we have the following commutative diagram
$$
\xymatrix{
 L(Y)\ar^{id}@{=}[rr]\ar[dr]
 &  &   L(Y).  \\
& L(Z) \ar[ur] }
$$
Since $L$ is full, we get the following diagram
$$
\xymatrix{
 Y\ar^{\beta}[rr]\ar_{x}[dr]
                &  &   Y     \\
                 & Z \ar_{y}[ur] }
$$
such that $L(\beta)=id=L(xy)$. Then we have a distinguished triangle
$Y\stackrel{xy}\ra Y\ra V\ra \Omega_{F}^{-1}(Y)$ in
$\mathscr{A}/\mathcal{P}(F)$. Thus, we get a distinguished triangle
$L(Y)\stackrel{id}\ra L(Y)\ra L(V)\ra L(Y)[1]$ in
$\DFb{\mathscr{A}}/\Kb{\mathcal{P}(F)}$ Consequently, $L(V)=0$.
Therefore, we have $V=0$ by our claim. It follows that $xy$ is an
isomorphism and $x$ is split monomorphism. Then $\alpha=0$ and $L$
is a faithful functor.

Step3. $L$ is dense. Let $X^{\bullet}$ be in
$\DFb{\mathscr{A}}/\Kb{\mathcal {P}(F)}$. Then $X^{\bullet}$ is in
$\DFb{\mathscr{A}}$. Take an $F$-projective resolution
$P_{X}^{\bullet}$ of $X^{\bullet}$ with
$P_{X}^{\bullet}\in\KFb{\mathcal {P}(F)}$. Without loss of
generality, we assume that $\cpx{P_{X}}$ has the following form:
$$\cdots\ra P_{X}^{-r-2}\ra P_{X}^{-r-1}\ra P_{X}^{-r} \stackrel{d^{-r}}\ra
\cdots \ra P_{X}^{0}\ra 0 \ra\cdots$$ such that
$H^{i}(Q,P_{X}^{\bullet})=0$ for $i<-r$ with any $Q\in \mathcal
{P}(F)$. Then we deduce that $$\cdots\ra P_{X}^{-r-2}\ra
Im(d_X^{-r-2})\ra 0 \ra\cdots $$ is an $F$-acyclic complex.
Therefore, the complex $\cpx{P_{X}}$ is isomorphic in
$\DFb{\mathscr{A}}$ to a complex of the form
$$ 0\ra Im(d_X^{-r-2})\ra
P_X^{-r-1}\stackrel{d^{-r-1}}\ra P_{X}^{-r}\ra\cdots\ra
P_{X}^{0}\ra0\ra\cdots.$$ From the following commutative diagram
$$\xymatrix{
\cdots \ar[r] & 0 \ar[r]\ar[d]& P_{X}^{-r-1}\ar[r]\ar@{=}[d] &
P_X^{-r}\ar[r]\ar@{=}[d]  &
\cdots \ar[r]\ar@{=}[d] & P_X^{0}\ar[r]\ar@{=}[d] & 0 \ar[r]\ar@{=}[d] &\cdots\\
\cdots \ar[r] & P_X^{-r-2}\ar[r]\ar@{=}[d]  & P_X^{-r-1}\ar[r]\ar[d]
& P_X^{-r}\ar[r]\ar[d] &
\cdots \ar[r]\ar[d] &P_X^{0}\ar[r]\ar[d]  & 0 \ar[r]\ar[d]&\cdots\\
\cdots \ar[r]&P_X^{-r-2}\ar[r]         & 0\ar[r] & 0 \ar[r]&
\cdots\ar[r] & 0\ar[r] & 0\ar[r]& \cdots,}$$ we conclude that
$Im(d_X^{-r-2})[r+2]\simeq P_{X}^{\bullet}\simeq X^{\bullet}$ in
$\DFb{\mathscr{A}}/\Kb{\mathcal {P}(F)}$. Take a minimal
$F$-injective resolution of $Im(d_X^{-r-2})$:
$$
0\ra Im(d_X^{-r-2})\ra I^{0}\ra I^{1}\ra\cdots\ra I^{r+1}\ra M\ra0
$$
such that $M\simeq \Omega_{F}^{-r-2}(Imd_X^{-r-2})$, where $I^{i}$
are $F$-injectives for $0\leq i\leq r+1$. Then
$Im(d_X^{-r-2})[r+2]\simeq M$ in $\DFb{\mathscr{A}}/\Kb{\mathcal
{P}(F)}$. Therefore, $L(M)\simeq X^{\bullet}$. This shows that $L$
is dense. $\square$

\noindent{\bf Remark.} In \cite{Gr}, Grime studied the special
closed subfunctor $F$ of $\Ext^1(-,-)$. We generalize his main
result to general closed subfunctor $F$.

\noindent{\bf Example.} The Artin algebra $\Lambda$ is said to be an
$F$-Frobenuis algebra if $\mathcal {P}(F)=\mathcal {I}(F)$. Set
$\mathscr{A}=\Lambda$-mod. Then we have a triangle equivalence
$\Lambda\text{-}$mod$/\mathcal {P}(F)\simeq
\DFb{\Lambda}/\Kb{\mathcal {P}(F)}$ by Theorem 5.1.

As a corollary of Theorem 5.1, we re-obtain Rickard's result.
\begin{Koro}$\rm\cite[Theorem\,2.1]{Ri2}$
Let $\Lambda$ be a self-injective Artin algebra. Then we have the
following triangle equivalence
$$\Db{\Lambda\text{-}mod}/\Kb{_{\Lambda}\mathcal {P}}\simeq
\smod\Lambda.$$
\end{Koro}


\section{Relative derived equivalences for Artin algebras}\label{sectionexchange}
In this section, suppose $\Lambda$ is an Artin $R$-algebras. Let $F$
be a non-zero subfunctor of $\Ext_{\Lambda}^{1}(-,-)$. Let
$\Lambda$-mod be the category of finitely generated left
$\Lambda$-modules. Denote by $\mathcal{P}(F)$ the subcategory of
$\Lambda$-mod consisting of all $F$-projective $\Lambda$-modules.
The subcategory $_{\Lambda}\mathcal{P}$ of finitely generated
projective $\Lambda$-modules is contained in $\mathcal{P}(F)$.
Assume that $F$ has enough projectives and injectives, such that
there exists $G\in \Lambda$-mod such that $\add G=\mathcal {P}(F)$.
Let $\DFb{\Lambda}$ be the relative bounded derived category of
$\Lambda$-mod.

We give our main result in this section, a Morita theory for
relative derived categories.

\begin{Theo}

 Let $\Lambda$ and $\Gamma$ be Artin algebras. The following conditions
are equivalent.

$(1)$ $\DFb{\Lambda}$ and $\Db{\Gamma}$ are equivalent as
triangulated categories.

$(2)$ $\Kb{\mathcal {P}(F)}$ and $\Kb{_{\Gamma}\mathcal {P}}$ are
equivalent as triangulated categories.

$(3)$ $\Kf{\mathcal {P}(F)}$ and $\Kf{_{\Gamma}\mathcal {P}}$ are
equivalent as triangulated categories.

$(4)$ $\Gamma$ is isomorphic to $\End(T^{\bullet})$, where
$T^{\bullet}$ is an object of
 $\Kb{\mathcal {P}(F)}$ satisfying

   \qquad $(i)$ $\Hom(T^{\bullet},T^{\bullet}[i])= 0$ for i $\neq$
   0.

   \qquad $(ii)$ $\add(T^{\bullet})$, the category of direct summands of
    finite direct sums of copies of $T^{\bullet}$, generates $\Kb{\mathcal {P}(F)}$ as
a triangulated category.
\end{Theo}

{\textbf{Proof}.} There exists a $\Lambda$-module $G$ such that
$\add G=\mathcal {P}(F)$. Let $\Sigma=\End_{\Lambda}(G)$. Since
$$
\Hom_{\Lambda}(G,-): \add G \ra _{\Sigma}\mathcal {P}
$$
is an equivalence, it follows that $\KFb{\mathcal {P}(F)}$ and
$\Kbz{_{\Sigma}\mathcal {P}}$ are equivalent as triangulated
categories by the functor $\Hom_{\Lambda}(G,-)$, where
$\KFb{\mathcal {P}(F)}=\{\cpx{X}\in \Kf{\mathcal {P}(F)}| \exists
\;n\in\mathbb{Z}\text{\; such that \;} H^i(\Hom_A(G,\cpx{X}))=0
\text{\; for\;} i\leq n \}.$ We want to show that $\Db{\Sigma}$ and
$\Db{\Gamma}$ are triangle equivalent. It suffices to prove that
$\Hom_{\Lambda}(G,T^{\bullet})$ is an ordinary tilting complex over
$\Sigma$.

(i) $\Hom_{\Lambda}(G,T^{\bullet})\in \Kb{_{\Sigma}\mathcal {P}}$.

(ii) $\Hom((G,T^{\bullet}),(G,T^{\bullet})[i]) \simeq
\Hom(T,T[i])=0$ $for i\neq 0$.

(iii) $\add(\Hom(G,T^{\bullet}))$ generates $\Kb{_{\Sigma}\mathcal
{P}}$ as triangulated category.

Finally,
$$
\Hom((G,T^{\bullet}),(G,T^{\bullet}))\simeq
\Hom(G\otimes_{\Sigma}(G,T^{\bullet}),T^{\bullet})\simeq
\Hom(T^{\bullet},T^{\bullet})=\Gamma.
$$
Then by \cite[Theorem 6.4]{Ri1}, $\Db{\Sigma}$ and $\Db{\Gamma}$ are
equivalent as triangulated categories. $\square$

\noindent{\bf Remark.} If $F=\Ext^{1}(-,-)$, then we get Rickard's
result \cite[Corollary 8.3]{Ri1}.

\begin{Def} If $\Lambda$ and $\Gamma$ satisfy the
equivalent conditions of Theorem 6.1, then $\Lambda$ and $\Gamma$
are said to be relatively derived equivalent.
\end{Def}

We shall call an object $T^{\bullet}$ of $\Kb{\mathcal {P}(F)}$
satisfying conditions (i) and (ii) of  Theorem 6.1(4) an $F$(or a
relative)-tilting complex for $\Lambda$.

Next, we give an example of a relative tilting complex.

\noindent{\bf Example.} Let $k$ be a field and $\Lambda$ the algebra
given by the quiver
$$
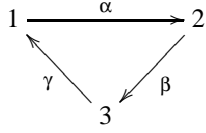
\xymatrix{
  1 \ar[rr]^{\alpha}
                &  &    2 \ar[dl]^{\beta}    \\
                & 3  \ar[ul]^{\gamma}               }
$$
and relations
$\alpha\beta\gamma=\beta\gamma\alpha\beta=\gamma\alpha\beta\gamma=0$.

Let $P_{1}, P_{2}, P_{3}$ be indecomposable projective
$\Lambda$-modules associated with vertices $1,2,3$, respectively.
Let $F=F_{\mathcal {P}(F)}$ be the subbifunctor with $\mathcal
{P}(F)=\add(P_{1}\oplus P_{2}\oplus P_{3}\oplus M)$, where
$M=\rad(P_{1})$. Let $T^{\bullet}_{1}$ be the complex $0\ra Q
\stackrel{f}\ra M\oplus P_{1}\ra 0$, where $f$ is an $\add
(P_{2}\oplus P_{3})$-approximation of $(M\oplus P_{1})$ and
$Q\in\add (P_{2}\oplus P_{3})$ is in degree $-1$. Let
$T^{\bullet}_{2}: 0\ra P_{2}\oplus P_{3}\ra0$ be a stalk complex
concentrated in degree $-1$. In this relative theory,
$T^{\bullet}=T^{\bullet}_{1}\oplus T^{\bullet}_{2}$ is an
$F$-tilting complex over $\Lambda$. This has been showed in
\cite[Proposition 1.2]{HK1} and \cite[Theorem 3.5]{HX1}. There is a
distinguished triangle $$Q[-1]\ra(M\oplus P_{1})[1]\ra\cpx{T}_{1}\ra
Q.$$ Then $G=P_{1}\oplus P_{2}\oplus P_{3}\oplus M\in\add\cpx{T}$.
By the proof of Theorem 6.1, it follows that $\Hom(G,\cpx{T})$ is a
tilting complex for $\End(G)$.

As for the theory of derived equivalence, it is useful to study the
Grothendieck groups in order to determined the number of
indecomposable non-isomorphic summands which a tilting complex has.
We introduce the following group. Let
$\mathbb{Z}(\Lambda\text{-}$mod) denote the free abelian group with
the isomorphism classes $[X]$ of $\Lambda$-modules $X$. Define the
$F$-stable Grothendieck group of $\Lambda$-mod,
$F$-$K_{0}(\Lambda$-mod), to be
$\mathbb{Z}(\Lambda$-mod)$/F$-$R(\Lambda$-mod), where
$F$-$R(\Lambda$-mod) is the subgroup of $\mathbb{Z}(\Lambda$-mod)
generated by the elements $[X]+[Z]-[Y]$ whenever $0\ra X\ra Y\ra
Z\ra 0$ is an $F$-exact sequence of $\Lambda$-modules.

\begin{Prop}Let $\Lambda$ and $\Gamma$ be relatively
derived equivalent. Then
$$
F\text{-}K_{0}(\Lambda\text{-}mod)\simeq K_{0}(\Gamma\text{-}mod).
$$
\end{Prop}

\textbf{Proof.} If $\Lambda$ and $\Gamma$ are relatively derived
equivalent, then $K_{0}(\DFb{\Lambda})\simeq K_{0}(\Db{\Gamma})$. As
is known, $K_{0}(\Gamma\text{-}mod)\simeq K_{0}(\Db{\Gamma})$. It
suffices to show that $F$-$K_{0}(\Lambda$-mod)$\simeq
K_{0}(\DFb{\Lambda})$.

Define a map $\alpha: F\text{-}K_{0}(\Lambda\text{-}mod)\ra
K_{0}(\DFb{\Lambda})$ by $\alpha([X])=[l(X)]$, where
$l:\Lambda\text{-}mod\ra \Db{\Gamma}$ is the inclusion functor.

Suppose that $X\ra Y\ra Z\ra X[1]$ is a distinguished triangle in
$\DFb{\Lambda}$ and $P\in \mathcal {P}(F)$. Applying the functor
$\Hom_{\DFb{\Lambda}}(P,-)$ to the above triangle, we have a exact
sequence
$$0\ra (P,X)\ra (P,Y)\ra (P,Z)\ra 0.
$$
It follows that $0\ra X\ra Y\ra Z\ra0$ is an $F$-exact sequence in
$\Lambda$-mod. Then we define $\beta: K_{0}(\DFb{\Lambda})\ra
F$-$K_{0}(\Lambda$-mod) by
$\beta([X^{\bullet}])=\Sigma_{i}(-1)^{i}[(X^{i})]$, where
$X^{\bullet} \in \DFb{\Gamma}$. It is easy to show that
$\alpha\beta=id_{K_{0}(\DFb{\Lambda})}$ and
$\beta\alpha=id_{F\text{-}K_{0}(\Lambda\text{-}mod)}$. $\square$

As a consequence of Proposition 6.3 we have the following
proposition.

\begin{Prop}
Let $T^{\bullet}$ be an $F$-tilting complex and
$\Gamma=\End(T^{\bullet})$. Then the following conditions are
equivalent.

$(1)$ The number of non-isomorphic indecomposable summands of
$T^{\bullet}$.

$(2)$ The number of non-isomorphic simple $\Gamma$-modules.

$(3)$ The number of non-isomorphic  indecomposable modules in
$\mathcal {P}(F)$.
\end{Prop}

By J.Rickard's criterion \cite[Proposition 6.2]{Ri1}, we give a
characterization of $\Kb{\mathcal {I}(F)}$ which is a full
triangulated subcategory of $\DFb{\Lambda}$.

\begin{Lem}\label{6.5} $\Kb{\mathcal {I}(F)}=\{X^{\bullet}\in
\DFb{\Lambda}\mid$ for any $Y^{\bullet}\in \DFb{\Lambda}$, there
exists $n_{0}\in\mathbb{Z}$, such that
$\Hom_{\DFb{\Lambda}}(Y^{\bullet},X^{\bullet}[i])=0, \forall i\geq
n_{0}\}$.
\end{Lem}

\textbf{ Proof.} If $X^{\bullet}\in \DFb{\Lambda}$ and $X^{\bullet}
\simeq I_{X}^{\bullet}$ with $I_{X}^{\bullet}\in \Kb{\mathcal
{I}(F)}$. Then
$$
\Hom_{\DFb{\Lambda}}(Y^{\bullet},X^{\bullet}[i])\simeq
\Hom_{\DFb{\Lambda}}(Y^{\bullet},I_{X}^{\bullet}[i])\simeq
\Hom_{\Kb{\Lambda}}(Y^{\bullet},I_{X}^{\bullet}[i])=0
$$
for all $i\geq n_{0}$ with some $n_{0}\in \mathbb{Z}$.

If $X^{\bullet}\in \DFb{\Lambda}$ and for any $Y^{\bullet}\in
\DFb{\Lambda}$, $\exists\;n_{0}\in\mathbb{Z}$, such that
$\Hom_{\DFb{\Lambda}}(Y^{\bullet},X^{\bullet}[i])=0$, for all $i\geq
n_{0}$. Since $\DFb{\Lambda}\simeq\Kbf{\mathcal {I}(F)}$, it follows
that $X^{\bullet}\simeq I_{X}^{\bullet}$ with $I_{X}^{\bullet}\in
K^{+,b}(\mathcal {I}(F))$ such that $H^{i}(P,I_{X}^{\bullet})=0$ for
$i\geq m, m\in\mathbb{Z}$. For $j\geq m$, we have the following
diagram
$$\xymatrix{
\cdots \ar[r]& I_{X}^{j-1}\ar[r]\ar@{=}[d]  &
I_{X}^{j}\ar[r]\ar@{=}[d]
& Im(d_X^{j})\ar[r]\ar[d]& 0\\
\cdots \ar[r]& I_{X}^{j-1}\ar[r] &  I_{X}^{j}\ar[r] &
I_{X}^{j+1}\ar[r]& \cdots,}$$ such that the complex
$I_{X}^{\bullet}$ is isomorphic in $\DFb{\Lambda}$ to a complex of
the form $ \cdots\ra I_{X}^{j-1}\ra I_{X}^{j}\ra Im(d_X^{j})\ra 0$.
We will prove that $Im(d_X^{j})$ is $F$-injective. If $Im(d_X^{j})$
is not $F$-injective, then there is $j=i-1\geq m$ and
$Im(d_X^{i-1})\neq 0$. Hence, the map $I_{X}^{i-1}\ra Im(d_X^{i-1})$
is not zero. Set $Y=\oplus_{j\geq m} Im(d_X^{j})$. Then the map
$I_{X}^{i-1}\ra Y=Im(d_X^{i-1})\oplus_{j\geq m, j\neq i-1}
Im(d_X^{j})$ is not zero. Consequently,
$\Hom_{\DFb{\Lambda}}(Y,X^{\bullet}[i])\neq0$ for $i\geq m+1$. This
leads to a contradiction with our assumption. $\square$

\section{Relative derived equivalences and relative homological dimensions}
In this section,we study the relationships between relative
dimensions and relative derived equivalences.

Let $\Lambda$ be an Artin algebra. Let us recall some notions of
relative homological dimensions. The $F$-projective dimension of a
$\Lambda$-module $X$, denoted by $\pd_{F}X$, is defined to be the
smallest $m$, such that there is an $F$-projective resolution
$$
0\ra P^{-m} \ra \cdots \ra P^{-1} \ra P^{0} \ra X \ra 0
$$
with $P^{-i} \in \mathcal {P}(F)$ for $i\geq 0$. If such $m$ does
not exist, then define $\pd_{F}X=\infty$. The $F$-global dimension
of $\Lambda$ is defined as
$$
\gldim_{F}(\Lambda)= sup\{\pd_{F}(X)\mid X\in \Lambda\text{-mod}\}.
$$
The $F$-finitistic dimension of $\Lambda$ is defined as follows
$$
\fd_{F}(\Lambda)= sup\{\pd_{F}X\mid\pd_{F}(X) <\infty, X\in
\Lambda\text{-mod} \}.
$$
The {\bf relative finitistic dimension} $rel.\fd_{F}(\Lambda)$ is
defined to be $sup_{F}\fd_{F}(\Lambda)$, where $F$ ranges over all
subfunctors $F=F_{\mathscr{C}}$ for all subcategories of finite type
of $\Lambda$-mod contains $\mathcal {P}(\Lambda)$.

We begin with the characterization of $F$-projective dimension of a
$\Lambda$-module.

\begin{Lem}The following conditions are equivalent.

$(i)$ $\pd_{F}(X)\leq m$.

$(ii)$ $\Ext^{i}_{F}(X,Y)=0$ for all modules $Y$ and for all $i\geq
m+1$.

$(iii)$ $\Ext^{m+1}_{F}(X,Y)=0$ for all modules $Y$.
\end{Lem}

Recall that $\Lambda$ and $\Gamma$ are relatively derived
equivalent, if there is a triangle equivalence $\DFb{\Lambda}\simeq
\Db{\Gamma}$, where $\Gamma\simeq\End(\cpx{T})$, $\cpx{T}$ is a
relative tilting complex for $\Lambda$. Clearly, $T^{\bullet}$ is a
bounded complex. Let $n\geq 0$ be an integer. Without loss of
generality, we may assume that $T^{\bullet}$ is a radical complex of
the form:
$$\cdots\rightarrow0\rightarrow T^{-n}\rightarrow
T^{-n+1}\rightarrow\cdots \rightarrow T^{-1}\rightarrow
T^{0}\rightarrow 0 \rightarrow\cdots.$$ The integer $n$ is called
the term length of $\cpx{T}$, denote by $t(T^{\bullet})=n$.

Let $H: \Db{\Gamma}\ra \DFb{\Lambda}$ be a relative derived
equivalence such that $H(\Gamma)=T^{\bullet}$. Then we have the
following fact.

\begin{Lem}Let $n\geq 0$ be an integer and $T^{\bullet}$ as above.
 Let $L: \DFb{\Lambda}\ra \Db{\Gamma}$ be a
quasi-inverse of $H$. Then there exists a radical tilting complex
$\bar{T}^{\bullet} \in \Kb{_{\Gamma}\mathcal {P}}$ such that
$L(G)=\bar{T}^{\bullet}$ and $\bar{T}^{i}=0$ for $i>n$ or $i<0$.
\end{Lem}

\textbf{Proof.} Considering the homological group of $\bar{T}$, we
have
$$
H^{i}(\bar{T}^{\bullet})\simeq
\Hom_{\Kb{\Gamma}}(\Gamma,\bar{T}^{\bullet}[i])\simeq
\Hom_{\Db{\Gamma}}(\Gamma,\bar{T}^{\bullet}[i])\simeq
\Hom_{\DFb{\Lambda}}(T^{\bullet},G[i]).
$$
By Proposition \ref{3.7}, we have
$\Hom_{\DFb{\Lambda}}(T^{\bullet},G[i])\simeq
\Hom_{\Kb{\Lambda}}(T^{\bullet},G[i])$. Then
$H^{i}(\bar{T}^{\bullet})=0$ for $i>n$ or $i<0$ by transferring
shifts. Since $\bar{T}^{\bullet}\in \Kb{_{\Gamma}\mathcal {P}}$, it
follows that $\bar{T}^{\bullet}$ has the following form
$$
\cdots\to 0\to \bar{T}^{-r}\to \bar{T}^{-r+1}\to \cdots \to
\bar{T}^{-1}\stackrel{d^{-1}}\to \bar{T}^{0} \stackrel{d^{0}}\to
\bar{T}^{1}\to \bar{T}^{s}\to0 \to\cdots
$$
where $r,s \in \mathbb{N}$. If $H^{i}(\bar{T}^{\bullet})=0$ for
$i>n$, then the following sequence $$0\to Im(d_T^{n})\to
\bar{T}^{n+1}\to \cdots \to \bar{T}^{s}\to0 \to\cdots$$ is a split
exact sequence. Consequently, $Im(d_T^{n})$ is a projective
$\Gamma$-module. Then the complex $\bar{T}^{\bullet}$ is isomorphic
in $\Db{\Gamma}$ to a complex of the form
$$0\to\bar{T}^{-r}\to \bar{T}^{-r+1}\to \cdots
\to \bar{T}^{-1}\to \bar{T}^{0}\to \cdots \to \bar{T}^{n-1}\to
Ker(d_T^{n})\to 0,
$$where $Ker(d_T^{n})$ is a projective $\Gamma$-module.
On the other hand, we have
$$
H^{i}(\Hom_{\Gamma}(\bar{T}^{\bullet},\Gamma))\simeq
\Hom_{\Db{\Gamma}}(\bar{T}^{\bullet}, \Gamma[i])\simeq
\Hom_{\DFb{\Lambda}}(G, T^{\bullet}[i])\simeq \Hom_{\Kb{\Lambda}}(G,
T^{\bullet}[i])=0.
$$for $i>0$.
It follows that $$0\to Im(d_T^{0},\Gamma)\to
(\bar{T}^{-1},\Gamma)\to (\bar{T}^{-2},\Gamma)\to\cdots$$ is split
exact in $\Kb{\Gamma^{op}}$, where
$Im(d_T^{0},\Gamma)=(Im(d_T^{-1}),\Gamma)$. We conclude that
$Ker(d_T^{0},\Gamma)=(Coker(d_T^{-1}),\Gamma)$ is a projective
$\Gamma$-module. Consequently, $Coker(d_T^{-1})$ and $Im(d_T^{-1})$
are projective $\Gamma$-modules. Therefore, the complex
$\bar{T}^{\bullet}$ is isomorphic a complex of the following form
$$
\cdots \to 0 \to Coker(d_T^{-1})\to \bar{T}^{1}\to \cdots \to
\bar{T}^{n-1}\to Ker(d_T^{n})\to 0 \quad\text{in}\quad
\Kb{_{\Gamma}\mathcal {P}}.
$$
This completes the proof. $\square$

The main result of this section can be stated as follows.

\begin{Theo}\label{7.3} Let $L:\DFb{\Lambda}\ra \Db{\Gamma}$ be a
relative derived equivalence. Suppose $T^{\bullet}$ is the relative
tilting complex associated to $L$. Then we have

$(1)$ $\gldim_{F}(\Lambda)-t(T^{\bullet})\leq \gldim(\Gamma)\leq
\gldim_{F}(\Lambda)+t(T^{\bullet})+2$.

$(2)$ $\fd_{F}(\Lambda)-t(T^{\bullet})\leq \fd(\Gamma)\leq
\fd_{F}(\Lambda)+t(T^{\bullet})+2$.
\end{Theo}

The following lemmas are useful in our proofs.

\begin{Lem}\label{7.4} Let $Y$ be a $\Gamma$-module. Then
$H^{i}(H(Y))=0$ for $i<-n$ or $i>0$, and in $\DFb{\Lambda}$, $H(Y)$
is isomorphic to a complex $\cpx{T_{Y}}$of the following form
$$
\cdots\ra 0\ra Im(d_{T_{Y}}^{-n-2})\ra T_{Y}^{-n-1}\ra \cdots \ra
T_{Y}^{-1}\ra Ker(d_{T_{Y}}^{0})\ra 0\ra \cdots,
$$
where $T_{Y}^{i}$ and $Ker(d_Y^{0})$ are $F$-projective
$\Lambda$-modules.
\end{Lem}

\textbf{Proof.}  We have the following isomorphisms
$$
H^{i}(H(Y))\simeq \Hom_{\DFb{\Lambda}}(\Lambda,H(Y)[i])\simeq
\Hom_{\Db{\Gamma}}(\bar{T}_{\Lambda}^{\bullet},Y[i]) \simeq
\Hom_{\K{\Gamma}}(\bar{T}_{\Lambda}^{\bullet},Y[i]),
$$
where $L(\Lambda)\simeq \bar{T}_{\Lambda}^{\bullet}$ has the same
form as $\bar{T}^{\bullet}$. Then $H^{i}(H(Y))=0$ for $i>0$ or
$i<-n$. Assume that $H(Y)$ is isomorphic to a bounded above
$F$-projective complex:
$$
\cdots\ra T_{Y}^{-n-1}\ra T_{Y}^{-n}\ra \cdots\ra T_{Y}^{-1}\ra
T_{Y}^{0}\ra\cdots\ra T_{Y}^{s}\ra0
$$
Then we get two acyclic complexes
$$
\cdots\ra T_{Y}^{-n-1}\ra Im(d_{T_{Y}}^{-n-1})\ra 0\ra \cdots
$$
and
$$
0\ra Im(d_{T_{Y}}^{0})\ra T_{Y}^{1}\ra \cdots.
$$
For any $P\in \mathcal {P}(F)$, we have
$$
H^{i}(P,H(Y))\simeq \Hom_{\DFb{\Lambda}}(P,H(Y)[i])\simeq
\Hom_{\Db{\Gamma}}(\bar{T}_{P}{^{\bullet}},Y[i]) \simeq
H^{i}(\Hom_{\Gamma}(\bar{T}_{P}{^{\bullet}},Y)) =0
$$
for $i>0$ or $i<-n$, where $L(P)\simeq \bar{T}_{P}{^{\bullet}}$ has
the same form as $\bar{T}^{\bullet}$. Therefore, we get the
following acyclic complexes:
$$
\cdots\ra \Hom(P,T_{Y}^{-n-1})\ra Im\Hom(P,d_{T_{Y}}^{-n-1})\ra 0
$$
and
$$
0\ra Im(\Hom(P,d_{T_{Y}}^{0}))\ra \Hom(P,T_{Y}^{1})\ra \cdots.
$$
We conclude that there exist two $F$-acyclic complexes:
$$
\cdots\ra T_{Y}^{-n-2}\ra Im(d_{T_{Y}}^{-n-2})\ra 0\ra \cdots,
$$
and
$$
0\ra Im(d_{T_{Y}}^{0})\ra T_{Y}^{1}\ra \cdots.
$$
Then the complex $H(Y)$ is isomorphic in $\DFb{\Lambda}$ to a
complex $\cpx{T_{Y}}$ of the following form:
$$
\cdots\ra 0\ra Im(d_{T_{Y}}^{-n-2})\ra T_{Y}^{-n-1}\ra \cdots \ra
T_{Y}^{-1}\ra Ker(d_{T_{Y}}^{0})\ra 0\ra \cdots,
$$
where $Ker(d_{T_{Y}}^{0})$ and $T_{Y}^{-i}$($i=0,\cdots, n-1$) are
$F$-projective $\Lambda$-modules. $\square$ \vspace{0.3cm}

\begin{Lem}\label{7.5} Let $X$ be a $\Lambda$-module. Then
$H^{i}(L(X))=0$ for $i>n$ or $i<0$, and in $\Db{\Gamma}$, $L(X)$ is
isomorphic to a complex $\cpx{\bar{T}_{X}}$ of the following form:
$$
\cdots\ra 0\ra Coker(d_{T_{X}}^{-1})\ra \bar{T}_{X}^{1}\ra \cdots
\ra \bar{T}_{X}^{n-1}\ra Ker(d_{T_{X}}^{n})\ra 0\ra \cdots
$$
where $\bar{T}_{X}^{i}$ ($1\leq i\leq n-1$) and $Kerd^{n}$ are
projective $\Gamma$-modules.
\end{Lem}

\begin{Lem}$\rm\cite[Lemma\,3.7]{P}$\label{7.6} Let $m,t,d \in \mathbb{N}$, $X^{\bullet}, Y^{\bullet} \in
\Kb{\Lambda\text{-}mod}$. Assume that $X^{i}=0$ for $i<m$, $Y^{j}=0$
for $j>t$, and $\Ext_{F}^{l}(X^{i},Y^{j})=0$ for all $i,j\in
\mathbb{N}$ and $l\geq d$. Then
$\Hom_{\DFb{\Lambda}}(X^{\bullet},Y^{\bullet}[l])=0$ for $l\geq
d+t-m$.
\end{Lem}

\textbf{Proof.} The proof is similar to the one in \cite[Lemma
1.6]{Ka}. The detailed proofs will appear in \cite{P}.
$\square$\vspace{0.3cm}

We now have all the ingredients to complete the proof of our main
theorem.

\textbf{Proof of Theorem 7.3:}

Assume that $\fd_{F}(\Lambda)<\infty$. Let $N$ be arbitrary
$\Gamma$-module, and let $M$ be a $\Gamma$-module with
$\pd_{\Gamma}(M)<\infty$. Let $P_{M}^{\bullet}\to M$ be a minimal
projective resolution of $M$. We have the following isomorphisms
$$
\Ext_{\Gamma}^{i}(M,N)\simeq \Hom_{\Db{\Gamma}}(M,N[i]) \simeq
\Hom_{\Db{\Gamma}}(P_{M}^{\bullet},N[i]) \simeq
\Hom_{\DFb{\Lambda}}(H(P_{M}^{\bullet}), H(N)[i]).
$$
By Lemma \ref{7.4}, we see that $H(N)$ is isomorphic in
$\DFb{\Lambda}$ to a complex $\cpx{T_{N}}$ of the form:
$$
\cdots\ra 0\ra Im(d_{T_{N}}^{-n-2})\ra T_{N}^{-n-1}\ra \cdots \ra
T_{N}^{-1}\ra Ker(d_{T_{N}}^{0})\ra 0\ra \cdots.
$$
By Theorem 6.1, it follows that $H(P_{M}^{\bullet})\simeq
T_{M}^{\bullet}\in \Kb{\mathcal {P}(F)}$ and
$H^{i}(H(P_{M}^{\bullet}))=0$ for $i<-n$ or $i>0$. Then we see that
$$
0\to T_{M}^{-r}\to \cdots \to Ker(d_{T_{M}}^{-n-2})\to 0
$$
and
$$
0\to Im(d_{T_{M}}^{0})\to \cdots \to T_{M}^{s}\to 0
$$
are $F$-split acyclic complexes, where $Im(d_{T_{M}}^{0})$ is an
$F$-projective $\Lambda$-module. Therefore, $T_{M}^{\bullet}$ is
isomorphic to a complex of the following form
$$0\to Im(d_{T_{M}}^{-n-2})\to T_{M}^{-n}\to \cdots \to
T_{M}^{-1} \to Ker(d_{T_{M}}^{0})\ra 0\quad\text{in}\quad
\DFb{\Lambda},$$ where $Ker(d_{T_{M}}^{0})$ and $T_{M}^{i}$ ($-n\leq
i\leq-1$) are $F$-projective $\Lambda$-modules. From the above
argument, we deduce that $\pd_{F}(Im(d_{T_{M}}^{-n-2}))<\infty $.
Consequently, $\pd_{F}(Im(d_{T_{M}}^{-n-2}))\leq \fd_{F}(\Lambda)$.
By Lemma \ref{7.6}, we conclude that
$$\Ext_{\Gamma}^{i}(M,N)\simeq\Hom_{\DFb{\Lambda}}(H(P_{M}^{\bullet}),H(N)[i])=0\quad\text{for}
\quad i\geq \fd_{F}(\Lambda)+1+n+2.
$$ This implies that $\fd(\Gamma)\leq
\fd_{F}(\Lambda)+l(T^{\bullet})+2$.

Assume that $\fd(\Gamma)<\infty$. Let $Y$ be any $\Lambda$-module,
and let $X$ be a $\Lambda$-module with $\pd_{F}(X)<\infty$. Let
$P_{X}^{\bullet}\to X$ be a minimal $F$-projective resolution of
$X$. We have the following isomorphisms
$$
\Ext_{F}^{i}(X,Y)\simeq \Hom_{\DFb{\Lambda}}(X,Y[i]) \simeq
\Hom_{\DFb{\Lambda}}(P_{X}^{\bullet},Y[i]) \simeq
\Hom_{\Db{\Gamma}}(L(P_{X}^{\bullet}), L(Y)[i]).
$$
By Lemma \ref{7.5}, we see that $L(Y)$ is isomorphic in
$\Db{\Gamma}$ to a complex of the following form:
$$
\cdots\ra 0\ra Coker(d_{T_{Y}}^{-1})\ra \bar{T}_{Y}^{1}\ra \cdots
\ra \bar{T}_{Y}^{n-1}\ra Ker(d_{T_{Y}}^{n})\ra 0\ra \cdots.
$$
From Theorem 6.1, it follows that $H(P_{X}^{\bullet})\simeq
\bar{T}_{X}^{\bullet}\in \Kb{_{\Gamma}\mathcal {P}}$ and
$H^{i}(L(P_{X}^{\bullet}))=0$ for $i<0$ or $i>n$. Therefore,
$\bar{T}_{X}^{\bullet}$ is isomorphic to a complex of the form
$$
\cdots\ra 0\ra Coker(d_{T_{Y}}^{-1})\ra \bar{T}_{X}^{1}\ra \cdots
\ra \bar{T}_{X}^{n-1}\ra Ker(d_{T_{Y}}^{n})\ra 0\ra
\cdots\quad\text{in}\quad \Db{\Gamma},
$$ where $Ker(d_{T_{Y}}^{n})$, $T_{X}^{i}$ ($1\leq i\leq n$) are projective
$\Gamma$-modules. We see that $\pd(Coker(d_{T_{Y}}^{-1}))<\infty$.
Consequently, $\pd(Coker(d_{T_{Y}}^{-1}))\leq \fd(\Gamma)$. By Lemma
\ref{7.6}, we conclude that
$$\Ext_{F}^{i}(X,Y)\simeq\Hom_{\Db{\Gamma}}(L(P_{X}^{\bullet}),L(Y)[i])=0\quad\text{for}\quad i\geq
\fd(\Gamma)+1+n.$$ This implies that $\fd_{F}(\Lambda)\leq
\fd(\Gamma)+l(T^{\bullet})$. $\square$

\noindent{\bf Remark.} From Theorem 7.3, we see that if the relative
finitistic dimension of $\Lambda$ is finite if and only if the
finitistic dimension of $\Gamma$ is finite. Auslander and Solberg
\cite[Proposition 6.1]{ASo2} proved that the standard finitistic
dimension conjecture is true for all Artin algebras if and only if
the relative finitistic dimension conjecture is true. Then we get
another approach to the finitistic dimension conjecture.

\noindent{\bf Example.} The standard global dimension
$\gldim(\Lambda)$ can be infinite, but $\gldim(\Gamma)$ is finite. A
less trivial example is the following. Let $\Lambda$ be the path
algebra of the following quiver
$$
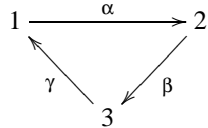
\xymatrix{
  1 \ar[rr]^{\alpha}
                &  &    2 \ar[dl]^{\beta}    \\
                & 3  \ar[ul]^{\gamma}               }
$$
and relations
$\alpha\beta\gamma=\beta\gamma\alpha=\gamma\alpha\beta=0$. Let
$P_{1}=\begin{smallmatrix}1\\2\\3\end{smallmatrix}$,
$P_{2}=\begin{smallmatrix}2\\3\\1\end{smallmatrix}$ and
$P_{3}=\begin{smallmatrix}3\\1\\2\end{smallmatrix}$, which are all
the indecomposable projective $\Lambda$-modules. Then let
$M_{i}=P_{i}/soc(P_{i})$ and $S_{i}=P_{i}/\rad(P_{i})$ for
$i=1,2,3$, where $soc(P_i)$ is the socle of $P_i$ for each $i$. Let
$\mathcal {P}(F)=\{P_{1}, P_{2}, P_{3}, S_2, S_3, M_{2}\}$. Then we
have $\mathcal {I}(F)=\D Tr({\mathcal {P}(F)})\cup
_{\Lambda}\mathcal {I}=\{P_{1}, P_{2}, P_{3}, S_{3},S_{1}, M_3\}$,
where $\D Tr$ is the Auslander-Reiten translation of $\Lambda$. It
is easy to find $F$-exact sequences for the modules $S_1$, $M_1$ and
$M_3$ to show that $\gldim_{F}(\Lambda)\leq1$. By Theorem 7.3 we
have $\gldim(\Gamma)\leq \gldim_{F}(\Lambda)+l(T^{\bullet})+2$ for
any $F$-tilting complex $\cpx{T}$, where
$\End(\cpx{T})\simeq\Gamma$. In particular, for the $F$-tilting
module
$$\cpx{T}=0\ra P_{1}\oplus P_{2}\oplus P_{3}\oplus S_2\oplus S_3\oplus
M_{2}\ra0,$$ where $P_{1}\oplus P_{2}\oplus P_{3}\oplus S_2\oplus
S_3\oplus M_{2}$ is in degree zero, we get $\gldim(\Gamma)\leq3$,
but $\gldim(\Lambda)<\infty$.

Recall that an Artin algebra $\Lambda$ is called an $F$-Gorenstein
algebra if $\mathcal {I}^{\infty}(F)=\mathcal {P}^{\infty}(F)$,
where $\mathcal {P}^{\infty}(F)=\{X\mid\pd_{F}(X)<\infty \}$,
$\mathcal {I}^{\infty}(F)=\{Y\mid\id_{F}(Y)<\infty \}$.

\begin{Lem}$\rm\cite[Lemma\,3.2]{ASo4}$ An Artin algebra $\Lambda$ is an
$F$-Gorenstein algebra if and only if both the relative injective
dimension of all modules in $\mathcal {P}(F)$ and the relative
projective dimension of all modules in $\mathcal {I}(F)$ are
bounded.
\end{Lem}

It is well known that Gorensteinness is invariant under derived
equivalence. It is not immediately obvious that the relation between
$F$-Gorensteinness and relative derived equivalences, but this will
follow from the following propositions.
\begin{Prop}\label{7.8}
 Let $L:\DFb{\Lambda}\ra \Db{\Gamma}$ be a relative derived
 equivalence. Suppose that $T^{\bullet}$ is the relative tilting complex
  associated to $L$. Then
$\id_{F}(\mathcal {P}(F))-t(T^{\bullet})\leq
\id(_{\Gamma}\Gamma)\leq \id_{F}(\mathcal {P}(F))+t(T^{\bullet})+2$.
\end{Prop}

\textbf{Proof.} Assume that $\id_{F}(\mathcal {P}(F))<\infty$. Set
$r=max\{\id(P)\mid P\in \mathcal {P}(F) \}$. By homological formula,
we have
$$
\Ext_{\Gamma}^{i}(X,\Gamma)\simeq\Hom_{D^{b}(\Gamma)}(X,\Gamma[i])
\simeq\Hom_{\DFb{\Lambda}}(H(X),T^{\bullet}[i]).
$$
By Lemma \ref{7.4}, $H(X)$ is isomorphic in $\DFb{\Lambda}$ to a
complex of the form:
$$
\cdots\ra 0\ra Im(d_{T_{X}}^{-n-2})\ra T_{X}^{-n-1}\ra \cdots \ra
T_{X}^{-1}\ra Ker(d_{T_{X}}^{0})\ra 0\ra \cdots.
$$
Moreover, $T^{\bullet}\in \Kb{\mathcal {P}(F)}$ and
$\id_{F}(\mathcal {P}(F))<\infty$. By Lemma 7.1, it follows that
$\Ext_{F}^{i}(-,P)=0$ for $i\geq r+1$. By Lemma \ref{7.6}, we deduce
that $\Ext_{\Gamma}^{i}(X,\Gamma)=0 $ for $i\geq r+1+n+2$. Then
$\id(_{\Gamma}\Gamma)\leq n+2+\id_{F}(\mathcal {P}(F))$.

Suppose that $\id(_{\Gamma}\Gamma)<\infty$. By homological formula,
we have
$$
\Ext_{F}^{i}(Y,P)\simeq\Hom_{\DFb{\Lambda}}(Y,P[i])\simeq\Hom_{\Db{\Gamma}}(L(Y),L(P)[i]).
$$
By Lemma \ref{7.5}, $L(Y)$ is isomorphic in $\Db{\Gamma}$ to a
complex of the form:
$$
\cdots\ra 0\ra Coker(d_{T_{Y}}^{-1})\ra \bar{T}_{Y}^{0}\ra \cdots
\ra \bar{T}_{Y}^{n-1}\ra Ker(d_{T_{Y}}^{n})\ra 0\ra \cdots
$$
By Lemma \ref{7.6}, we conclude that $\Ext_{F}^{i}(Y,P)=0$ for
$i\geq \id(_{\Gamma}\Gamma)+1+n$. Consequently, $ \id_{F}(\mathcal
{P}(F))\leq \id(_{\Gamma}\Gamma)+n$. This completes the proof.
$\square$ \vspace{0.3cm}

\begin{Prop}\label{7.9}
Let $\Lambda$ and $\Gamma$ be relatively derived equivalent Artin
$R$-algebras. Then $\Lambda$ is $F$-Gorenstein if and only if
$\Gamma$ is Gorenstein.
\end{Prop}

\textbf{Proof.} Assume that $\Lambda$ is an $F$-Gorenstein algebra.
Then we have $\id_{F}(\mathcal {P}(F))<\infty$ and $\pd_{F}(\mathcal
{I}(F))<\infty$. Set $r=max\{\pd_{F}(J)\mid J\in \mathcal {I}(F)\}$.
For any injective $\Gamma$-module $I$, we have
$$
\Ext_{\Gamma}^{i}(I,X)\simeq\Hom_{\Db{\Gamma}}(I,X[i])\simeq\Hom_{\DFb{\Lambda}}(H(I),H(X)[i]).
$$
By Lemmas \ref{6.5} and \ref{7.4}, we conclude that $H(I)\in
\Kb{\mathcal {I}(F)}$ and $H(I)$ is isomorphic in $\DFb{\Lambda}$ to
a complex $\cpx{J}$ of the following form:
$$
\cdots\ra 0\ra Im(d_{J}^{-n-2})\ra J^{-n-1}\ra \cdots \ra J^{-1}\ra
J^{0}\ra \cdots\ra J^{s}\ra0
$$
such that $Im(d_{J}^{-n-2})$, $J^{-i}$ are $F$-injective
$\Lambda$-modules for $-n-1\leq i\leq s$. By Lemma \ref{7.6}, we see
that $\Ext_{\Gamma}^{i}(I,X)=0 $ for $i\geq r+1+n+2$. This implies
that $\pd(I)\leq n+2+\pd_{F}(\mathcal {I}(F))$. Therefore,
$\id(\Gamma_\Gamma)<\infty$. By Proposition \ref{7.8}, it follows
that if $\id_{F}(\mathcal {P}(F))<\infty$, then
$\id(_{\Gamma}\Gamma)< \infty$. Therefore, $\Gamma$ is a Gorenstein
algebra.

If $\Gamma$ is a Gorenstein algebra, then
$\id(_{\Gamma}\Gamma)=\id(\Gamma_{\Gamma})=s<\infty$. In the same
way, we deduce that $\Lambda$ is an $F$-Gorenstein algebra. This
completes the proof. $\square$ \vspace{0.3cm}

In particular, when $F=\Ext_{\Lambda}^{1}(-,-)$, by the proof of
Theorem \ref{7.3}, we have the following.

\begin{Koro} Let $\Lambda$ and $\Gamma$ be derived equivalent.
Suppose $T^{\bullet}$ is a tilting complex for $\Lambda$. Then

$(1)$ $\gldim(\Lambda)-l(T^{\bullet})\leq \gldim(\Gamma)\leq
\gldim(\Lambda)+l(T^{\bullet})$.

$(2)$ $\fd(\Lambda)-l(T^{\bullet})\leq \fd(\Gamma)\leq
\fd(\Lambda)+l(T^{\bullet})$.

$(3)$ $\id(\Lambda)-l(T^{\bullet})\leq \id(\Gamma)\leq
\id(\Lambda)+l(T^{\bullet})$.
\end{Koro}

\noindent{\bf Acknowledgements.} The author gratefully thanks
Professor Changchang Xi for his many helpful suggestions and would like to
acknowledge National Natural Science Foundation
of China (No. $11201022$) and the Fundamental Research Funds for the Central Universities (2013JBM096, 2013RC027). The author
also would like to thank the referee for his/her helpful comments
that improve the paper.

\bigskip

\begin{thebibliography}{99}
{\small
\bibitem{ARS}{{\sc  Auslander, M., Reiten, I., Smal{\O}, S. O.}:
Representation Theory of Artin Algebras. Cambridge University
Press, 1995.}
\bibitem{ASo1}{{\sc Auslander, M., Solberg, {\O}.:}
Relative homology and representation theory I, Relative homology and
homologically finite subcategories. \emph{Comm. Algebra} \textbf{21}(1993),
2995-3031.}
\bibitem{ASo2}{{\sc Auslander, M., Solberg, {\O}.:}
Relaive homology and representation theory II, Relative cotilting
theory. \emph{Comm. Algebra}  \textbf{21}(1993), 3033-3079.}

\bibitem{ASo3}{{\sc Auslander, M., Solberg, {\O}.:}
Relative homology and representation theory III, Cotilting modules
and Wedderburn correspondence. \emph{Comm. Algebra} \textbf{21}(1993),
3081-3097.}

\bibitem{ASo4}{{\sc Auslander, M., Solberg, {\O}.:}
Gorenstein algebras and algebras with dominant dimension at least
$2$. \emph{Comm. Algebra}  \textbf{21}(1993), 3897-3934.}

\bibitem{ASo5}{{\sc Auslander, M., Solberg, {\O}.:}
Relative homology.\, Finite-dimensional algebras and related topics
(Ottawa, ON, 1992), 347--359, NATO Adv. Sci. Inst. Ser. C Math.
Phys. Sci., 424, Kluwer Acad. Publ., Dordrecht, 1994.}

\bibitem{Bu1}{{\sc Buan, A.:}
Closed subbifunctors of the extension bifunctor. \emph{J. Algebra}
\textbf{244}(2001), 407-428.}

\bibitem{Bu2}{{\sc Buan, A.:}
Subcategories of the derived category and cotilting complexes. \emph{Colloq. Math.} \textbf{88}(2001), 1-11.}

\bibitem{Buch}{{\sc Buchweitz,  R. O.:}
Maximal Cohen-Macaulay modules and Tate-cohomology over Gorenstein
rings, Hamburg, p. 155 (1987). (unpublished manuscript)}

\bibitem{Ch}{{\sc Chen, X. W.:}
Gorenstein homological algebra of Artin algebras, Postdoctoral
Report, USTC, 2010.}

\bibitem{CYZ}{{\sc Chen, X. W., Ye, Y., Zhang, P.:}
Algebras of dereived dimension of zero. \emph{Comm. Algebra} \textbf{36}(2008),
1-10.}

\bibitem{CZ}{{\sc Chen, X. W.,  Zhang, P.:}
Quotient triangulated categories. \emph{Manuscripta Math.} \textbf{123}(2007),
167-183.}
\bibitem{CBS}{{\sc Cline, E., Parshall, B., Scott L.:}
Derived categories and Morita theory. J. Algebra 104(1986),
397-409.}
\bibitem{DRSS}{{\sc P. Draxler}, {\sc I. Reiten}, {\sc S. O. Smal{\O}} and {\sc {\O}. Solberg},
Exact categories and vector space categories. \emph{Trans. Amer.
Math. Soc.} \textbf{351}(1999), 647-682.}

\bibitem{DS}{{\sc Dugger, D., Shipley, B.:}
K-theory and derived equivalences. \emph{Duke Math. J.} \textbf{124}(2004),
587-617.}

\bibitem{EJ}{{\sc Enochs, E. E., Jenda, O. M. G.:}
Relative Homological Algebra. De Gruyter Expositions in
Math. 30, Walter de Gruyter-Berlin-New York, 2000.}

\bibitem{GZ}{{\sc Gao, N., Zhang, P.:}
Gorenstein derived categories. \emph{J. Algebra} \textbf{323}(2010), 2041-2057}.

\bibitem{GM}{{\sc Gelfand,  S. I., Manin, Y.:}
Methods of Homological Algebra. Springer Monographs in
Math. (2nd ed.), Springer-Verlag, 2002.}

\bibitem{Ge}{{\sc Generalov, A. I.:}
Relative homological algebra. Cohomology of categories, posets and
coalgebras, \, in "Handbook of Algebra"\, V.1,\, 611-638.\,
North-Holland, Amsterdam, 1996.}

\bibitem{Gr}{{\sc Grime, M.:}
Adjoint functors and triangulated categories. \emph{Comm. Algebra}
\textbf{36}(2008), 3589-3607.}

\bibitem{Ha1}{{\sc Happel, D.:}
Triangulated Categories in the Representation Theory of Finite
Dimensional Algebras. Cambridge University Press, Cambridge. 1988.}

\bibitem{Ha2}{{\sc Happel, D.:}
On the derived category of finite-dimensional algebra. \emph{Comment.
Math. Helv.} \textbf{62}(1987), 339-389.}

\bibitem{Ha3}{{\sc Happel, D.:}
On Gorenstein algebras.\, In: Representation theory of finite groups
and finitedimensional algebras \, (Proc. Conf. at Bielefeld, 1991),
389-404, \, Progress in Math., vol. 95, Birkh$\ddot{a}$user, Basel,
1991.}


\bibitem{Har}{{\sc Harshorne, R.:}
Duality and Residue. Lecture Notes in Math. 20, Springer,
Berlin, 1966.}

\bibitem{Hoch}{{\sc Hochschild, G.:}
Relative homological algebra. \emph{Trans. Amer. Math. Soc.} \textbf{82}(1956),
246-269.}

\bibitem{HK1}{{\sc Hoshino, M., Kato, Y.:}
Tilting complexes defined by idempotents. \emph{Comm. Algebra} \textbf{30}(2002),
83-100.}
\bibitem{HK2}{{\sc Hoshino, M., Kato, Y.:}
Tilting complexes associated with a sequence of idempotents. \emph{J.
Algebra} \textbf{183}(2003), 105-124.}

\bibitem{HX1}{{\sc Hu, W., Xi, C. C.:}
$D$-split sequences and derived equivalences.
\emph{Adv. Math.} \textbf{227}(2011), 292-318.}

\bibitem{HX2}{{\sc Hu, W., Xi, C. C.:}
Derived equivalences and stable equivalences of Morita type, I.
\emph{Nagoya Math. J.} \textbf{200}(2010), 107-152.}

\bibitem{HX3}{{\sc Hu, W., Xi, C. C.:}
Derived equivalences for $\Phi$-Auslander-Yoneda algebras. \emph{Trans. Amer. Math. Soc.} \textbf{365}(2013), 5681-5711.}


\bibitem{KS}{{\sc Kashiwara, M., Schapira, P.:}
Sheaves on Manifolds. Grundlehren der mathematischen
Wissenschaften 292\, Berlin: Springer-Verlag, 1990.}

\bibitem{Ka}{{\sc Kato, Y.:}
On derived equivalent coherent rings. \emph{Comm. Algebra} \textbf{30}(2002),
4437-4454.}

\bibitem{Ke1}{{\sc Keller, B.:}
Deriving DG categories. \emph{Ann. Sci. \'{e}cole Norm. Sup.} \textbf{27}
(1994), 63-102.}


\bibitem{KZ}{{\sc Konig, S.,Zimmermann, A.:}
Derived Equivalences for Group Rings. Lecture notes in
math. 1685, Berlin: Springer-Verlag, 1998}

\bibitem{Kr2}{{\sc Krause, H.:}
Localization theory for triangulated categories. \, To appear in the
proceedings of the "Workshop on Triangulated Categories" in Leeds,
2006.}

\bibitem{Mi2}{{\sc Miyashita, Y.:}
Tilting modules associated with a series of idempotent ideals. \emph{J.
Algebra} \textbf{238} (2001), 485¨C501.}

\bibitem{N}{{\sc Neeman, A.}
Triangulated Categories. Annals of Mathematics Studies 148,
Princeton University Press, \, Princeton and Oxford,\, 2001.}

\bibitem{P}{{\sc Pan, S. Y.:}
Relative quotient triangulated categories. Accepted by Algebra Colloquium.}

\bibitem{PX}{{\sc Pan, S. Y., Xi, C. C.:}
Finiteness of finitistic dimension is invariant under derived
equivalences. \emph{J. Algebra} \textbf{322}(2009), 21-24.}


\bibitem{Ri1}{{\sc Rickard, J.:}
Morita theory for derived categories. \emph{J. London Math. Soc.} \textbf{39}
(1989), 436-456.}

\bibitem{Ri2}{{\sc Rickard, J.:}
Derived categories and stable equivalence. \emph{J. Pure Appl. Algebra} \textbf{61}
(1989), 303-317.}

\bibitem{Ri3}{{\sc Rickard, J.:}
Derived equivalences as derived functors. \emph{J. London Math. Soc.}  \textbf{43}
(1991), 37-48.}

\bibitem{RoT}{{\sc Rotman, J.:}
An Introduction to Homological Algebra. Academic Press New
York, San Francisco, London, 1979.}


\bibitem{TT}{{\sc Thomason, R. W., Trobaugh, T. F.:}
Higher algebraic $K$-theory of schemes and of derived categories.
Then Grothendieck Festschrift (a collection of papers to honor
Grothendieck's 60th birthday) vol 3, Birkh\"{a}user, 1990, pp.
247-435.}


\bibitem{Ver}{{\sc Verdier, J.:}
Cat\'{e}gories d\'{e}riv\'{e}es, \'{e}tat 0. \, Lecture Notes in
Math. 569(1977), Springer, Berlin, 262-311.}

\bibitem{Xi1}{{\sc Xi, C. C.:}
The relative Auslander-Reiten theory of modules. Preprint, available
at : http://math.bnu.edu.cn/~ccxi/Papers/Articles/rart.pdf, 2002.}


\bibitem{W}{{\sc Weibel, C. A.:}
An Introduction to Homological Algebra. Cambridge Studies
in Advanced Math. 38,  Cambridge University Press, 1994.}

}
\end{thebibliography}
\end{document}